\documentclass[a4paper,16pt]{article}
\usepackage{amsthm}
\usepackage[sumlimits]{amsmath}
\usepackage{amsfonts}
\usepackage{color}
\usepackage{graphics,latexsym}
\allowdisplaybreaks

\def\MM#1{\boldsymbol{#1}}
\newtheorem{theorem}{Theorem}[section]

\newtheorem{proposition}[theorem]{Proposition}
\newtheorem{remark}[theorem]{Remark}
\DeclareMathOperator{\diff}{d}
\DeclareMathOperator{\Id}{Id}

\DeclareMathOperator{\ad}{ad}

\DeclareMathOperator{\Diff}{Diff}

\DeclareMathOperator{\Vol}{Vol}
\DeclareMathOperator{\VPM}{VPM}

\newcommand{\pp}[2]{\frac{\partial #1}{\partial #2}} 
\newcommand{\dd}[2]{\frac{\delta #1}{\delta #2}}
\newcommand{\dede}[2]{\frac{d #1}{d #2}}
\newcommand{\dt}[1]{\diff\!#1}

\newtheorem{corollary}[theorem]{Corollary}
\newtheorem{definition}[theorem]{Definition}
\newcommand{\bfi}{\bfseries\itshape}
\newcommand{\rem}[1]{}
\begin{document}

\title{Discrete momentum maps for lattice EPDiff} 
\author{\bigskip
  Colin J. Cotter${}^{1}$ and Darryl D. Holm${}^{1,2}$
  \\${}^{1}$Mathematics Department, \\
  Imperial College London, SW7 2AZ, UK \\{\footnotesize email:
    colin.cotter@imperial.ac.uk, d.holm@imperial.ac.uk} \bigskip
  \\${}^{2}$Computer and Computational Science,  \\
  Los Alamos National Laboratory, Los Alamos, NM 87545, USA
  \\{\footnotesize email: dholm@lanl.gov} }
\maketitle
\begin{abstract}
  We focus on the spatial discretization produced by the Variational
  Particle-Mesh (VPM) method for a prototype fluid equation the known
  as the {\it EPDiff equation}, which is short for {\it
    Euler-Poincar\'e equation associated with the diffeomorphism group
    \textup{(}of $\mathbb{R}^d$, or of a $d$-dimensional manifold
    $\Omega$\textup{)}}.  The EPDiff equation admits measure valued
  solutions, whose dynamics are determined by the momentum maps for
  the left and right actions of the diffeomorphisms on embedded
  subspaces of $\mathbb{R}^d$. The discrete VPM analogs of those
  dynamics are studied here. Our main results are: (\emph{i}) a
  variational formulation for the VPM method, expressed in terms of a
  constrained variational principle principle for the Lagrangian
  particles, whose velocities are restricted to a distribution
  $D_{\VPM}$ which is a finite-dimensional subspace of the Lie algebra
  of vector fields on $\Omega$; (\emph{ii}) a corresponding
  constrained variational principle on the fixed Eulerian grid which
  gives a discrete version of the Euler-Poincar\'e equation; and
  (\emph{iii}) discrete versions of the momentum maps for the left
  and right actions of
  diffeomorphisms on the space of solutions.

\end{abstract}
\tableofcontents

\section{Introduction}

\subsection{Transverse internal wave interactions}
\begin{figure}[htp]
\begin{center}
\scalebox{0.35}{\includegraphics{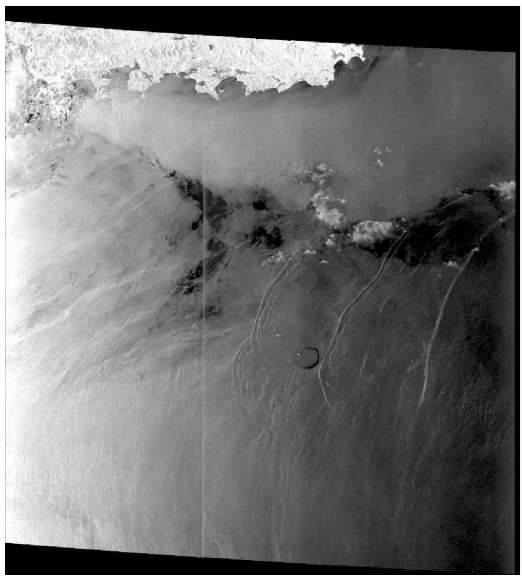}}
\end{center}
\caption{\label{dongsha}Image from the Space Shuttle of long,
  tidally-excited waves in the South China Sea near the Dong sha
  Atoll.  The waves are propagating from East to West, and are
  produced every tide (about 12 hours). The waves interact with the
  Atoll and then undergo nonlinear reconnections. Picture courtesy
  of A. Liu.}
\end{figure}

Synthetic Aperture Radar (SAR) observations from the Space Shuttle often show
nonlinear internal wave trains that propagate for many hundreds of
kilometers across large basins such as the South China Sea (SCS)
shown in Figure \ref{dongsha}.

These wave trains are characterized as \emph{Great Lines on the Sea}
in \cite{yoder94}. Both lines and spirals on the sea arise as flow
phenomena, rather than wave phenomena \emph{per se}. The flow
phenomenon detected in the the SAR imagery is associated with
nonlinear internal waves, whose crests may be as much as 200km long.
The amplitude of these internal waves results in about 150m of
deflection in the thermocline over a distance of about 1 km. Thus,
their aspect ratio satisfies the first criterion to be nonlinear
shallow water waves. Their amplitude is also considerably less than
the typical thickness of the thermocline, but it is not actually
infinitesimal compared to the thermocline thickness. The flow along
the crests of these waves also indicates they are not precisely the
same as usual shallow water waves.

The particular nonlinear internal waves found in the SCS are generated
by the tides flowing East to West through the Luzon Strait over
submerged ridges between Taiwan and the Phillipines. The SAR images in
Figure \ref{dongsha} show that the momentum of the tides flowing
Westward over these ridges concentrates into internal waves on the
thermocline that emerge into the SCS basin as thin wave fronts which
may extend in length for hundreds of kilometers (much larger than the
Straits in which they were created) and may propagate for thousands of
kilometers.  Perhaps because of the complex topography, the tides
flowing over the mouth of the Luzon Strait do not produce internal
waves propagating in both directions. The significant wave trains
propagate Westward.

\begin{figure}[htp]
\begin{center}
\scalebox{0.35}{\includegraphics{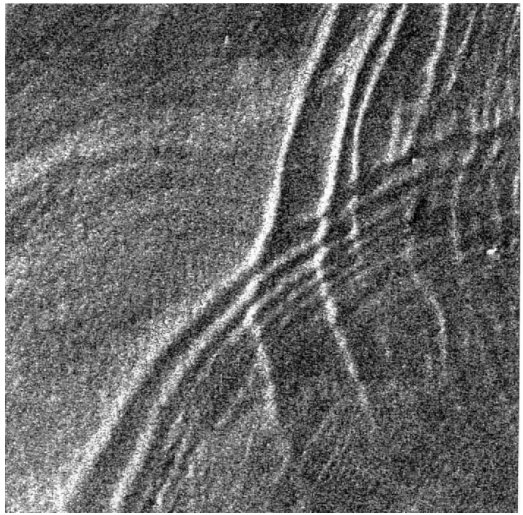}}
\end{center}
\caption{\label{zoom dongsha} Enlargement of part of figure
  \ref{dongsha} showing reconnecting long waves. Picture courtesy of A. Liu.}
\end{figure}

Propagating wave trains may intersect transversely with other wave
trains.  Sometimes these wave trains merely pass through each other as
linear waves.  However, in nonlinear wave encounters such as those
captured by SAR imaging of the region of the SCS West of Dong Sha
Island in Figure \ref{zoom dongsha}, two wave fronts may intersect
transversely, merge together and produce a single wave front.  This
merger of the wave fronts is the hallmark of a nonlinear process.
These particular wave interactions possess strong transverse dynamics
(flow along the crests) and momentum exchange in the direction of
propagation, which allow the wave fronts to merge and reconnect,
rather than merely passing through each other, as weaker waves do when
they intersect in an interference pattern.

Nonlinear internal wave interactions have been well studied in one
dimension, often by using the weakly nonlinear Boussinesq
approximation. These studies have usually resulted in a variant of the
Korteweg-de~Vries (KdV) equation, which has soliton solutions that
interact by exchange of momentum in unidirectional elastic collisions
(Whitham~1967).  However, the complex wave front interactions shown in
in Figure \ref{zoom dongsha} are plainly at least two-dimensional.  We
shall pursue the qualitative description of these higher-dimensional
wave interactions by using a simple two-dimensional model equation
called EPDiff.%
\footnote{EPDiff is the ``Euler-Poincare equation on the
  diffeomorphisms''.}  EPDiff may be derived in one dimension from the
asymptotic expansion for shallow water wave motion of the Euler
equations for the unidirectional flow of an incompressible fluid with
a free surface moving under gravity. In one dimension the result is
the Camassa-Holm (CH) equation, which arises at quadratic order in
this expansion. That is, CH is one order of accuracy in the asymptotic
expansion beyond KdV, which arises at linear order. Just as for KdV,
the CH equation is completely integrable; so CH also has soliton
solutions that interact by elastic collisions in one dimension.
Moreover, in the limit of zero linear dispersion, the CH solitons
develop a sharp peak at which their profile has a jump in derivative
that forms a sharp peak.  In this limit, the CH solitons are called
``peakons.'' The CH peakons are weak solutions, in the sense that
their momentum is concentrated on delta functions that move with the
velocity of the fluid flow.

In its zero-dispersion limit, CH has a geometric property that allows
it to be immediately generalized to higher dimensions, in which it is
called EPDiff. The term ``EPDiff'' distinguishes CH, which is a
one-dimensional shallow water wave equation with physical wave
dispersion, from its dispersionless limit which belongs to a larger
class of equations. This larger class of equations -- the
Euler-Poincar\'e (EP) equations \cite{HoMaRa1998} -- describes
geodesic motion with respect to any metric defining a norm on the
vector space of the Lie algebra of a Lie group. In the geometric
theory of fluid mechanics, the fluid velocity belongs to the tangent
space of the group of smooth invertible maps, called
``diffeomorphisms'' (or diffeos, for short).  The Euler-Poincare
equation on the diffeomorphisms is called EPDiff. EPDiff is a larger
class of equations than CH also because it is defined for geodesic
motion on the diffeos with respect to any metric, not just for the
$H^1$ norm of the velocity, which appears as the kinetic energy norm
in the derivation of CH. (The gradient part of the $H^1$ norm for CH
corresponds to the vertically averaged kinetic energy associated with
vertical motion.)  Thus, among the EPDiff equations, the
dispersionless limit of CH is one-dimensional EPDiff($H^1$). In one
dimension, the momentum of the EPDiff($H^1$) peakons is concentrated
at points moving along with the flow; but in higher dimensions, their
momentum is distributed on embedded subspaces moving with the flow. In
particular, EPDiff($H^1$) in two dimensions has singular solutions
whose momentum is distributed along curves in the plane. As solutions
of the two-dimensional version of a unidirectional shallow water wave
equation in its limit of zero linear dispersion, these moving curves
in the plane evolving under the dynamics of EPDiff($H^1$) are
prototypes for studying the interactions of the Great Lines on the
Sea.

To jump ahead, the singular (or, weak) solutions of the EPDiff
equation that emerge in finite time from any confined smooth initial
conditions and are supported on embedded subspaces moving with the
flow velocity, just as seen in the Great Lines on the Sea captured in
Figure \ref{dongsha}. We developed a numerical method for simulating
the singular solutions of EPDiff in the framework of its geometric
definition, which is natural for the Variational Particle Mesh (VPM)
method. Our numerical results using VPM show that
\begin{itemize}
\item
Singular solutions for EPDiff may be simulated by VPM as curve-segments
moving with the 2D flow velocity that possess no internal degrees of freedom.
\item
In collisions between any two of these curve-segment
solutions for EPDiff, the momentum of the one that overtakes from behind is
imparted to the one ahead. Thus, overtaking collisions between two
finite-length wave packets are elastic.
\item
The transverse collision of two curve-segment solutions for EPDiff may
result in merger (or, reconnection) of the  curve segments due to a
combination of exchange of momentum between the wave trains and flow along
their wave crests. In two dimensions, the reconnection or merger of singular
wave fronts under numerical EPDiff dynamics using VPM is evident in Figure
\ref{rearend}.
\end{itemize}

\paragraph{Plan of the paper} In this paper we introduce the VPM
method for EPDiff, and discuss some of the properties that arise from
the variational structure, in the following sections:
\begin{itemize}
\item The particle-mesh calculus is set out in
section \ref{setup-sec}.  
\item We give a variational principle associated
with the method in section \ref{HPDEP-sec}.
\item Section \ref{DEP-sec} shows that the Eulerian grid quantities
  satisfy an approximation the the EPDiff equation in Euler-Poincar\'e
  form.
\item Section \ref{L-Action-sec} defines a left action of $D_{\VPM}$
on $\Omega$ and a provides the corresponding momentum map.
\item Sections \ref{R-Action-sec} defines a right action in an
  extended space which can be interpreted as a discrete form of
  relabelling of Lagrangian particles. The Hamiltonian for the
  continuous time evolution of discretised EPDiff solutions is
  invariant under the action and so from Noether's theorem we obtain a
  conserved momentum.
\item Section \ref{kelvin} shows how this conserved momentum can be
  interpreted as a discrete form of Kelvin's circulation theorem.
\item Section \ref{numerics} gives some numerical examples, as well
as convergence tests for the method.
\end{itemize}

\subsection{Theoretical development}
Much of the theoretical development in this paper is inspired by the
following theorem \cite{Ar1966}.
\begin{theorem}[Arnold (1966) \cite{Ar1966}]
The solutions of Euler's equations for the incompressible motion of an
ideal fluid describe coadjoint geodesic motion on the volume preserving
diffeomorphisms, with respect to the $L^2$ norm of the fluid velocity (the
kinetic energy). 
\end{theorem}
The Euler equations for incompressible motion of an ideal fluid may be
written in the material frame as 
\[
\mathcal{P}\Big(\frac{d\mathbf{u}}{dt}\Big)=0
\quad\hbox{along}\quad
\frac{d\mathbf{x}}{dt}=\mathbf{u}
\quad\hbox{with}\quad
\nabla\cdot\mathbf{u}=0
\,,\]
where $\mathcal{P}$ is the Leray projection onto the incompressible vector
fields. These equations may also be written in the spatial frame as 
\[
\mathcal{P}\Big(\partial_t\mathbf{u} +{\rm
  ad}^*_\mathbf{u}\mathbf{u}\Big)=0 \,,\] where ad$^*$ is the dual of
the ad-action among incompressible vector fields under the $L^2$
pairing. That is, ad$^*$is defined by $\langle \MM{\mu}\,,\,{\rm
  ad}_\mathbf{u}\mathbf{w}\rangle = -\langle{\rm
  ad}^*_\mathbf{u}\MM{\mu}\,,\,\mathbf{w}\rangle$. Here $\ad_{\MM{u}}
\MM{w}=[\MM{u},\MM{w}]$ is the Lie-algebra commutator between vector
fields $\mathbf{u},\,\mathbf{w}$, and
$\langle\,\cdot\,,\,\cdot\,\rangle$ denotes the
$L^2$ pairing between such vector fields and
one-form densities such as   $\MM{\mu}$.

\paragraph{EPDiff}$\quad$\\

The EPDiff equation describes the corresponding coadjoint geodesic
motion on the {\it full} diffeomorphism group, allowing for
compressibility and an arbitrary norm, $\|\cdot\|$,
\[
\hbox{EPDiff is}\quad
\partial_t\MM{\mu}
+{\rm ad}^*_{\MM{u}}\MM{\mu}=0
\,,\quad\hbox{with}\quad
\MM{\mu}=\frac{\delta \ell}{\delta \MM{u}}
\,,\quad\hbox{where}\quad
\ell=\frac{1}{2}\|\MM{u}\|^2
\,.
\]
The momentum density $\MM{\mu}$ is a one-form density and the EPDiff
equation describes coadjoint dynamics under the action of the
corresponding velocity vector field. In EPDiff, ${\rm ad}^*_{\MM{u}}$ is the
coadjoint action of a vector field $\MM{u}$ acting on a one-form density
$\MM{\mu}=\delta\ell/\delta\MM{u}$ for a Lagrangian $\ell[\MM{u}]$ in
Hamilton's principle $\delta{S}=0$ for $S=\int\ell[\MM{u}]dt$.  In
components,
\[
\MM{\mu}=\mathbf{m}\cdot d\mathbf{x}\otimes dVol
\]
and EPDiff may be written as the invariance condition,
\[
\frac{d\MM{\mu}}{dt}=0
\quad\hbox{along}\quad
\frac{d\mathbf{x}}{dt}=\mathbf{u}=G*\mathbf{m}
\,,
\]
where $G*$ denotes convolution with the Green's function relating
the components of $\mathbf{m}$ and $\mathbf{u}$. In particular, for the
$H^1$ norm $\|\MM{u}\|^2\equiv\int |\mathbf{u}|^2
+\alpha^2|\nabla\mathbf{u}|^2\,dVol$, we have the component relation
\begin{equation}
\label{helmholtz}
\mathbf{m}=\mathbf{u}-\alpha^2\Delta\mathbf{u}
\,,
\end{equation}
and $G$ is the Green's function for the Helmholtz operator,
$Id-\alpha^2\Delta$, $\Delta$ is the Laplacian, and $\alpha$ is a
lengthscale. Thus, EPDiff for the $H^s$ norm with $s>0$ is an
integro-partial differential equation.

Originally derived \cite{HoMaRa1998} as an $n$-dimensional
generalisation of the Camassa-Holm equation for shallow-water dynamics
in one dimension \cite{CaHo1993}, EPDiff arises in several other
applications. For example, EPDiff for the $H^1$ norm is the
pressureless version of the Lagrangian-averaged Navier-Stokes-alpha
(LANS-alpha) model of turbulence \cite{FoHoTi2001}.  EPDiff for $H^1$
also emerges in the limit in which one ignores variations in height of
the Green-Nagdhi equation for shallow water dynamics
\cite{CaHoLe1996}. In one dimension, this is the dispersionless limit
of the Camassa-Holm equation \cite{CaHo1993}. In general, EPDiff is
the equation for coadjoint geodesic motion on the diffeomorphisms with
repect to any given norm on the Eulerian particle velocity (kinetic
energy).  Finally, EPDiff also describes the process of template
matching in computational anatomy \cite{MiTrYo2002}. In this
application, EPDiff has recently become a conduit for technology
transfer from soliton theory to computational anatomy
\cite{HoRaTrYo2004}. Thus, EPDiff turns out to be a prototype equation
for a number of applications.

The present article describes the underlying principles for using the
Variational Particle-Mesh (VPM) method in numerically integrating
EPDiff in the study of its nonlinear wave interactions.

\paragraph{Variational Particle-Mesh (VPM) method}
The Variational Particle-Mesh (VPM) method introduced in \cite{VPM}
produces Hamiltonian spatial discretizations of fluid equations which
may then be integrated in discrete time by using a variational
integrator. VPM may be regarded as a descendant of the Hamiltonian
Particle-Mesh method \cite{FrGoRe2002}, which is a Hamiltonian
discretization of the rotating shallow-water equations. The difference
is that HPM combines an Eulerian representation of the potential
energy (which gives rise to the pressure term) with a Lagrangian
representation of the kinetic energy, whilst VPM uses an Eulerian
representation of the entire Lagrangian. This means that the VPM
method is much more general than HPM and may be applied to many
different fluid PDEs (\emph{e.g.,} shallow-water, Green-Nagdhi,
incompressible Euler, \emph{etc.}). In this paper we focus on EPDiff,
which is an equation for fluid velocity only. Consequently, symmetries
of the discretised fluid velocity will be symmetries of the equations.
In future we will extend this work to include advected quantities such
as density, scalars \emph{etc.} Our ultimate aim is to use geometric
properties in constructing general numerical methods for PDEs
describing the continuum dynamics of fluids, complex fluids and
plasmas.

\rem{
which have
exact conservation of discrete circulation (which is the conservation
law arising from the particle-relabelling symmetry).
}

The conservative properties of variational integrators
are well understood \cite{LeMaOrWe2003}. In this article, we will discuss
preservation under VPM spatial discretization of the geometric
properties of the well-known EPDiff equation for coadjoint motion
under the diffeomorphisms \cite{HoMa2004},
\[
\partial_t\MM{\mu} + {\rm ad}^*_u\MM{\mu}=0
\,, \mbox{ with } \MM{\mu} = \dd{l}{\MM{u}}.
\]
In particular, we shall discuss discrete VPM analogs of the momentum
maps for the left and right actions of the diffeomorphisms on embedded
subspaces of $\mathbb{R}^n$ \cite{HoMa2004}.  The Lagrangian we shall
choose is the $H^1$ norm,
$\ell[\MM{u}]=\frac{1}{2}\|\MM{u}\|_{H^1}^2$, so the components of
velocity $\MM{u}$ and momentum density
$\MM{\mu}=\delta\ell/\delta{\MM{u}}$ will be related by the Helmholtz
operator, as in equation (\ref{helmholtz}). In this case, velocity
$\MM{\MM{u}}\in H^1$ implies that its dual momentum density
$\MM{\mu}\in H^{-1}$; so the solutions of EPDiff may be measure valued
in $\MM{\mu}$. That is, weak solutions of EPDiff are allowed in this
case, which are expressed in terms of delta functions supported on
embedded subspaces of $\mathbb{R}^n$ \cite{HoMa2004}. The left action
of the diffeomorphisms on these embedded subspaces of $\mathbb{R}^n$
generates the motion of spatially discrete EPDiff (lattice EPDiff),
while the right action is a symmetry and generates the conservation
law for circulation according to the Kelvin-Noether theorem
\cite{HoMaRa1998}. Thus, we seek the spatially discrete version of the
corresponding theorem for continuum solutions in \cite{HoMa2004}.  All
of these properties will then be preserved by an appropriate
variational time integrator.

\section{Particle-mesh calculus} \label{setup-sec} This section
describes the particle-mesh calculus that will be used in discretising
EPDiff. We shall describe its discretisation in space with continuous
time, and later we shall describe how to construct variational time
integrators to assemble a fully discrete space-time integration
scheme.

\paragraph{A finite dimensional subspace of $\mathfrak{X}(\Omega)$}
The infinite-dimensional space of smooth vector fields
$\mathfrak{X}(\Omega)$ generates the diffeomorphisms (smooth
invertible maps with smooth inverse) of the domain $\Omega$ onto
intself. To make a numerical algorithm that can be calculated on a
finite computer, we first need to choose a finite-dimensional subspace
$\mathfrak{X}_0$ of $\mathfrak{X}(\Omega)$ that will generate our
diffeomorphisms. We begin with a fixed grid consisting of $n_g$ points
in the domain $\Omega$ with vector coordinates
$\{\MM{x}_k\}_{k=1}^{n_g}\in\mathbb{R}^{d\times n_g}$ in $d$
dimensions. At each grid point $\MM{x}_k$ we shall associate a
velocity vector $\MM{u}_k$. The finite-dimensional space of possible
sets of velocity vectors $\MM{u}=(\MM{u}_1,\ldots,\MM{u}_{n_g})
\in\mathbb{R}^{d\times n_g}$ will then represent the required subspace
$\mathfrak{X}_0$. We call a set of values $\{\MM{u}_k\}_{k=1}^{n_g}$
the {\bfi grid representative} of the corresponding vector field.

To obtain the element of $\mathfrak{X}(\Omega)$ corresponding to 
$\MM{u}$, we use a set of basis functions with $\psi_k(\MM{x})$
representing a distribution centred around $\MM{x}_k$.  These basis
functions are taken to have compact support and to satisfy the {\bfi
Partition-of-Unity (PoU) property}
\[
\sum_{k=1}^{n_g}\psi_k(\MM{x})=1, \qquad \forall\ \MM{x}\in {\Omega}
\,.
\] 
The vector field $\MM{X}_{\MM{u}}$ is then defined as follows:
\begin{definition}
The vector field $\MM{X}_{\MM{u}}$ on $\Omega$ whose grid representative
is $\MM{u}$ takes the coordinate form
\[
\MM{X}_{\MM{u}}(\MM{x}) = \sum_k\MM{u}_k\psi_k(\MM{x})\cdot
\pp{}{\MM{x}}\,.
\]
\end{definition}

\begin{figure}[htp]
\begin{center}
\scalebox{0.7}{\includegraphics{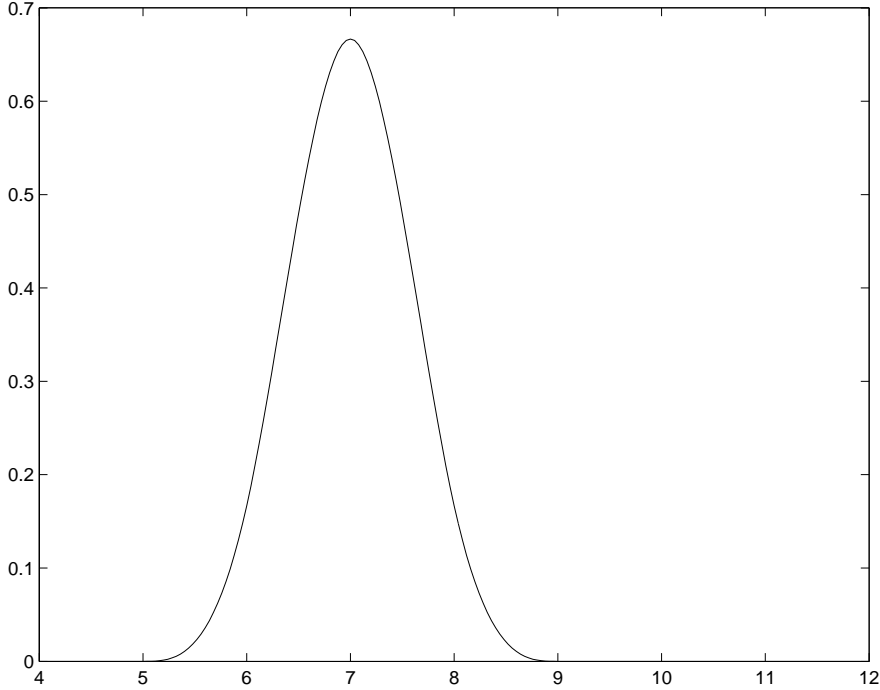}}
\end{center}
\caption{\label{blobs}Plot of a $B$-spline basis function $\psi_k$
  centred on $\MM{x}_k=7$ with a grid width of 1. These basis
  functions satisfy the partition-of-unity property.}
\end{figure}

A plot of a typical basis function in one dimension is given in figure
\ref{blobs}.

\begin{remark}
In general, these vector fields do not commute amongst themselves in the
Lie bracket, so they do not form a Lie subalgebra of
$\mathfrak{X}(\Omega)$. This will lead to a variational principle with
nonholonomic constraints. Also in general, the value of
$\MM{X}_{\MM{u}}(\MM{x}_k)$ is not exactly equal to $\MM{u}_k$, but is
convergent to it in the continuum limit. 
\end{remark}

\paragraph{Dynamics of a finite set of Lagrangian particles}
We shall proceed in describing our numerical method by introducing
 a finite set of $n_p$ {\bfi Lagrangian fluid particles}
$\{\MM{Q}_\beta\}_{\beta=1}^{n_p}$, whose velocities
$\{\MM{\dot{Q}}_\beta\}_{\beta=1}^{n_p}$ are entirely determined by the
grid velocity representation $\{\MM{u}_k\}_{k=1}^{n_g}$ \emph{via}
the vector field $\MM{X}_{\MM{u}}$ as follows: 
\begin{definition}
The {\bfi PoU vector field} $\MM{X}_{\MM{u};n_p}\in T\Omega^{n_p}$
associated with a velocity grid representative
$\MM{u}$ is defined as
\begin{equation}
\label{vector field def}
\MM{X}_{\MM{u};n_p}(\MM{Q}) = \sum_{\beta}\sum_k\MM{u}_k\psi_k(\MM{Q})
\cdot\pp{}{\MM{Q}}
\,,\quad
\MM{Q}\in(\MM{Q}_1,\ldots,\MM{Q}_{n_p})=\mathbb{R}^{d\times n_p}\,.
\end{equation}
\end{definition}
We shall constrain the dynamics of the particles so that a tangent
vector $\MM{\dot{Q}}$ may be represented as a PoU vector field evaluated
at the point $\MM{Q}$. That is, $(\MM{Q},\MM{\dot{Q}})$ lies in
a distribution $D^{\VPM}$ defined as follows:
\begin{definition}[The distribution $D^{\VPM}$]
Let $D^{\VPM}\subset T\Omega^{n_p}$ be the distribution defined by
\[
D^{\VPM} = \left\{
\left(\MM{\dot{Q}},\MM{Q}\right):
\MM{\dot{Q}}_\beta = \sum_k\MM{u}_k\psi_k(\MM{Q}_\beta)
\mbox{ for some }\MM{u}\in\mathbb{R}^{d\times n_g}
\mbox{ and }\forall\, k=1,\ldots,n_g
\right\}\,.
\]
\end{definition}
\begin{definition}\label{VPM-traj}
  A time series $\MM{Q}(t) = (\MM{Q}_1(t),\ldots,\MM{Q}_{n_p}(t))$
  with $(\MM{\dot{Q}}(t),\MM{Q}(t))\in D^{\VPM}\,\,\forall\, t_0\leq t
  \leq t_1$ is called a {\bfi VPM trajectory}. Each VPM trajectory
  defines a time series $\MM{u}_k(t)\in\mathbb{R}^{d\times n_g}\times
  [t_0,t_1]$ such that
  \begin{equation}
    \label{dot Q}
    \MM{\dot{Q}}_\beta(t) = \sum_k\MM{u}_k(t)\psi_k(\MM{Q}_\beta)
    \,,
  \end{equation}
  for $\beta = 1, \ldots,n_p$\,.  This is the {\bfi VPM tangent vector
    relation}, which we will enforce as a constraint for the
  variational principle resulting in the VPM method.
\end{definition}

\begin{remark}
  Given $(\MM{\dot{Q}},\MM{Q})$ one may invert equation (\ref{dot Q})
  for the grid velocity representation
  $\MM{u}=(\MM{u}_1,\ldots,\MM{u}_{n_g})\in\mathbb{R}^{d\times
    n_g}\times\mathbb{R}$, modulo the kernel of $\psi_k(\MM{Q}_\beta)$
  regarded as a matrix. (We shall see that this kernel does not affect
  the dynamics.) Later we shall write the Lagrangian as a function of
  $\MM{u}_k$ only and rely on this inversion to express the
  Euler-Lagrange equations for $\MM{Q}(t)$. We also note that a VPM
  trajectory $\MM{Q}(t)$ is specified entirely by
  $\{\MM{u}_k(t)\}_{k=1}^{n_g}$ and the initial condition $\MM{Q}(0)$.
  Changes of the initial conditions $\MM{Q}(0)$ for the VPM
  trajectories that leave invariant the grid velocity representation
  $\MM{u}(t)$ will provide the analog for VPM of ``particle
  relabelling'' in the continuum case.
\end{remark}

\paragraph{Gradient and divergence}

In this section we describe how the operations of gradient, divergence
and curl may be approximated using the particle-mesh discretisation. These
approximations apply the {\bfi two dual purposes of the basis functions
$\psi_k$}:
\begin{enumerate}
\item The $\psi_k$ interpolate functions from the grid to the particles.
\item The $\psi_k$ also construct densities on the grid from weights
stored on the particles.
\end{enumerate}
{\bf Notation:} Square brackets $[\,\cdot\,]^G$ and
$[\,\cdot\,]^P$ will denote these two  maps from particles to grid and
\emph{vice versa}. Superscripts distinguish whether the quantity is
evaluated on the grid or on the particles. That is, $[\,\cdot\,]^P$
indicates  mapping from grid to particles, and $[\,\cdot\,]^G$ indicates
mapping from particles to grid.
\begin{definition} Let $\{f_k\}_{k=1}^{n_g}$ be a scalar quantity 
stored at the Eulerian grid points. Then
\[
[f]_\beta^P = \sum_kf_k\psi_k(\MM{Q}_\beta),
\]
is an approximation of $f$ evaluated at the particle locations. 
Furthermore,
\[
[\nabla f]_\beta^P = \sum_kf_k\pp{}{\MM{Q}_\beta}\psi_k(\MM{Q}_\beta),
\]
is an approximation of the gradient of the scalar $f$ evaluated at the
particle locations.
\end{definition}
\begin{definition}
Let $\{g_\beta\}_{\beta=1}^{n_p}$ be a distribution
of values stored at particle locations. We construct a density
on the Eulerian grid as
\[
[g]_k^G = \sum_l(M^{-1})_{kl}\sum_\beta g_\beta\psi_l(\MM{Q}_\beta)\,.
\]
Furthermore, if the distribution is vector-valued $\MM{g}$ then
\[
[\nabla\cdot\MM{g}]_k^G =
-\sum_l(M^{-1})_{kl}\sum_\beta\MM{g}_\beta\cdot
\pp{}{\MM{Q}_\beta}\psi_l(\MM{Q}_\beta),
\]
is an approximation to the divergence of $\MM{g}$ on the grid.
\end{definition}
\paragraph{Discretised continuity equation} Given a set of constant
weights $\tilde{D}_\beta$ on the particles $\beta=1,\ldots,n_p$, 
to construct a density
\[
D_k = \sum_l(M^{-1})_{kl}\sum_\beta\tilde{D}_\beta\psi_k(\MM{Q}_\beta),
\]
one computes
\begin{eqnarray*}
  \dede{D_k}{t} & = & \sum_l(M^{-1})_{kl}\sum_\beta\tilde{D}_\beta
  \pp{\psi_l}{\MM{Q}_\beta}
  (\MM{Q}_\beta)\cdot\dot{\MM{Q}}_\beta, \\
  & = & \sum_l(M^{-1})_{kl}\sum_\beta\tilde{D}_\beta
  \pp{\psi_l}{\MM{Q}_\beta}\cdot\sum_k\MM{u}_k\psi_k(\MM{Q}_\beta),
\end{eqnarray*}
and so
\[
\dede{}{t}[D]^G=
-\,[\nabla\cdot\left([\MM{u}]^PD\right)]^G,
\]
so the corresponding grid representative $[D]^G$ satisfies a
discretised continuity equation.

\subsection*{Lagrangian for semi-discrete EPDiff}

Next we form the Lagrangian for semi-discrete EPDiff, as an
approximation to the continuous EPDiff Lagrangian
$$
L_C = \frac{1}{2}\int_{\Omega} (||\MM{u}||^2 + 
\alpha^2||\nabla\MM{u}||^2)\,\diff{\Vol},
$$
in which the constant $\alpha$ has dimensions of length.


\begin{definition}
  \label{lagrangian def}
  Let $\{N_k(\MM{x})\}_{k=1}^{n_g}$ be chosen as a finite element basis so
  that functions may be approximated in the form
  $$
  f(\MM{x}) = \sum_kN_k(\MM{x})f_k, \qquad f_k = f(\MM{x}_k)\,.
  $$
  (This basis need not be the same as that used in the
  partition-of-unity representation of velocity.) Define the matrix
  $H$, which approximates applying the Helmholtz operator and
  integrating, as
  $$ H_{kl} = \int_{\Omega} N_k(\MM{x})N_l(\MM{x}) + \alpha^2\nabla
  N_k(\MM{x}) \cdot\nabla N_l(\MM{x}) \diff{\Vol},
  $$ for some value of the constant $\alpha$. Then the Lagrangian for
  discrete EPDiff is expressed in this basis as
  \begin{equation}\label{normu-Lag}
    L(\MM{u}) = \frac{1}{2}\sum_{k,l}\MM{u}_k\MM{\cdot}
    H_{kl}\MM{u}_l \equiv \frac{1}{2}\MM{u}\cdot H\MM{u}
    \,.
  \end{equation}
\end{definition}
\begin{remark}
As in the continuous case, this Lagrangian is written entirely in
terms of the Eulerian velocity (in this case, the velocity grid
representation). In the continuum case, this form of the
Lagrangian admits Euler-Poincar\'e reduction (as Eulerian velocity is
invariant under the right-action of the diffeomorphism group
$\Diff(\Omega)$). This reduction results in the EP equation
\[
\MM{\MM{\mu}}_t+\ad^*_{\MM{u}}\MM{\MM{\mu}}=0
\,,
\] 
where $\MM{\MM{\mu}}=\delta L/\delta \MM{u}$ and ad$^*$ is the
\emph{dual} of the ad-action (Lie algebra commutator) of vector fields
on the domain. In the VPM discretisation of EPDiff, an analogous
equation will emerge, written on the Eulerian grid.
\end{remark}

\section{Variational principle for discrete EPDiff}\label{HPDEP-sec}
In this section we shall derive the equations for $\MM{Q}(t)$ from a
variational principle applied to the Lagrangian
(\ref{normu-Lag}) and required to satisfy the VPM tangent vector
constraint. Namely, the variational principle is constrained to
restrict the solutions so that
$\MM{\dot{Q}}\in D^{\VPM}_{\MM{Q}}$ (defined as the subspace
$\{\MM{\alpha}: (\MM{\alpha},\MM{Q})\in D^{\VPM}\}\subset
T_{\MM{Q}}\Omega^{n_p}$). This constraint on VPM trajectories
is the discrete analog of the {\bfi Lin constraints} in the
{\bfi Clebsch variational approach} to continuum ideal fluid
dynamics, as discussed for example in \cite{HoKu1983,HoMaRa1998}.
At the end of this section, we shall give a fully discrete
variational principle which produces the numerical scheme.

\subsection{Constrained action principle for semi-discrete EPDiff}
We begin by defining the {\bfi grid momentum} as follows. 
\begin{definition}[Grid momentum]
  \label{grid momentum}
  In the same finite element basis as for definition \ref{lagrangian
    def}, define the {\bfi integration matrix} $M$ (often called the
{\bfi mass matrix} in the finite element literature) as
  $$
  M_{kl} = \int_{\Omega}N_k(\MM{x})N_l(\MM{x})\diff{\Vol}.
  $$
  The {\bfi grid momentum} $\MM{m} = (\MM{m}_1,\ldots,\MM{m}_{n_g})
  \in \mathbb{R}^{d\times n_g}$ is then defined from the grid
  representative Lagrangian $L(\MM{u})$ in (\ref{normu-Lag}) \emph{via}
  its derivative
  \begin{equation}\label{gridmom-def}
  \sum_lM_{kl}\MM{m}_l = \pp{L}{\MM{u}_k}
  \,.
  \end{equation}
\end{definition}
\noindent
This expression for the grid momentum is an approximation to
$\delta{L}/\delta{\MM{u}}$ in the continous case.


\begin{definition}[Constrained action]
\label{implicit action principle}
The action for semi-discrete EPDiff is defined in terms of
three variables: the grid velocity $\MM{u}\in\mathbb{R}^{d\times
n_g}$; the particle positions $\MM{Q}_\beta\in\Omega^{n_p}$; and the
Lagrange multipliers $\MM{P}_\beta\in T^*_{\MM{Q}}\Omega^{n_p}$ which
will become the particle momenta on the Hamiltonian side. The action is 
given by
\begin{eqnarray}
\mathcal{A} = \int_0^T L(\MM{u}) + \sum_\beta\MM{P}_\beta\cdot
\left(
\MM{\dot{Q}}_\beta-\sum_k\MM{u}_k\psi_k(\MM{Q}_\beta)
\right)\diff{t}\,.
\label{VPM action}
\end{eqnarray}
This is the action for Lagrangian (\ref{normu-Lag})
when its particle velocities are required to satisfy the VPM tangent
vector constraint given in (\ref{dot Q}).
\end{definition}
\begin{proposition}\label{Prop-extreme}
The variables $(\MM{u},\MM{P},\MM{Q})$ which extremise the
constrained action $\mathcal{A}$ in (\ref{VPM action}) satisfy
\begin{eqnarray}
\MM{\dot{Q}}_\beta & = & \sum_k\MM{u}_k\psi_k(\MM{Q}_\beta)\,, 
\label{dotQ}\\
\label{dotP}
\MM{\dot{P}}_\beta & = & -\,\MM{P}_\beta\cdot\sum_k\MM{u}_k
\pp{\psi_k}{\MM{Q}}(\MM{Q}_\beta)\,, \\
\pp{L}{\MM{u}_k} & = & \sum_\beta\MM{P}_\beta\psi_k(\MM{Q}_\beta)
\,.
\label{gridmom-eqn}
\end{eqnarray}
\end{proposition}
\begin{proof}
After integration by parts, the first variation of $\mathcal{A}$ in
$(\MM{u},\MM{P},\MM{Q})$ is
\begin{eqnarray*}
\delta \mathcal{A} & = & \int_0^T \sum_k\left(
\pp{L}{\MM{u}_k} - \sum_\beta\MM{P}_\beta\psi_k(\MM{Q}_\beta)\right)
\cdot\delta\MM{u}_k \\
& & \quad +\sum_\beta\left(
\MM{\dot{Q}}_\beta-\sum_k\MM{u}_k\psi_k(\MM{Q}_\beta)\right)\cdot
\delta\MM{P}_\beta \\
& & \quad -\sum_\beta\left(
\MM{\dot{P}}_\beta+\MM{P}_\beta\cdot\sum_k\MM{u}_k\pp{\psi_k}{\MM{Q}}
(\MM{Q}_\beta)
\right)\cdot\delta\MM{Q}_\beta
\diff{t}\,,
\end{eqnarray*}
and the result follows by direct calculation.
\end{proof}
\begin{remark}{\bf [Left momentum map]}
\label{left mom-map}
Equation (\ref{gridmom-eqn}) in proposition \ref{Prop-extreme} bears a
great resemblance to the momentum map for left action of the
diffeomorphisms on embedded subspaces
\cite{HoMa2004} which describes the singular solutions of continuum EPDiff
equation. We will see later that equation (\ref{gridmom-eqn}) is the
discrete version of that momentum map.
\end{remark}

\begin{remark}\label{grid mom}{\bf [Grid momentum]}
The grid-momentum relation (\ref{gridmom-eqn}) allows one to obtain
$\MM{u}$ from $\partial{L}/\partial{\MM{\dot{Q}}}$ and $\MM{Q}$ by first
calculating $\MM{m}$, and then inverting the matrix $H_{kl}$ in 
\begin{equation}
\label{discrete elliptic}
\sum_lM_{kl}\MM{m}_l = \pp{L}{\MM{u}_k} = \sum_{l}H_{kl}\MM{u}_l
\,.
\end{equation}
This is the discrete analogue of the problem of solving for
$\MM{u}$ from $\MM{m}$ in the elliptic relation
$$
\MM{m} = \dd{L}{\MM{u}} = (1-\alpha^2\Delta)\MM{u}\,.
$$
for the continuous case \cite{HoMa2004}. Thus, the Lagrangian
(\ref{normu-Lag}) is hyper-regular on the grid. 
\end{remark}

\subsection{Legendre transform}
We now pass to the Hamiltonian side \emph{via} the Legendre
transform, a process summarised in the following proposition.
\begin{proposition}
  \label{ham side}
 The system of equations (\ref{dotQ}-\ref{gridmom-eqn}) is 
 canonically Hamiltonian with Hamiltonian function
$H$ given by
  \begin{equation}\label{Ham-mom}
  H = \frac{1}{2}\sum_{k,l}M_{kl}\MM{m}_l\cdot\MM{u}_k
  \,,
  \end{equation}
  where $\MM{u}_k=\sum_{m,n}H^{-1}_{km}M_{mn}\MM{m}_n$ is defined in
terms of $\sum_lM_{kl}\MM{m}_l$
  by inverting the matrix $H_{kl}$ in equation (\ref{discrete elliptic}),
and where $\MM{m}_k$ is obtained \emph{via} equations
(\ref{gridmom-def}) and (\ref{gridmom-eqn}).
\end{proposition}
\begin{proof}
  We obtain the Hamiltonian \emph{via} the Legendre transform
  $$
  H(\MM{P},\MM{Q}) = 
  \sum_\beta\MM{P}_\beta\cdot\sum_k\MM{u}_k\psi_k(\MM{Q}_\beta) -
  L(\MM{u}) \,,$$ 
  subject to equation (\ref{gridmom-eqn}).  Upon applying equation
  (\ref{gridmom-eqn}) the phase space action sum may be written as
  \begin{eqnarray*}
    \sum_\beta\MM{P}_\beta\cdot\sum_k\MM{u}_k\psi_k(\MM{Q}_\beta) &=& 
    \sum_k\MM{u}_k\cdot\sum_\beta\MM{P}_\beta\psi_k(\MM{Q}_\beta) 
    \\
    &=& \sum_{k,l}\MM{u}_k\cdot M_{kl}\MM{m}_l
    \,,
  \end{eqnarray*}
  after switching the orders of summation. The Lagrangian
  (\ref{normu-Lag}) may also be
  rewritten as
  $$
  L = \frac{1}{2}\MM{u}\cdot A\MM{u} = \frac{1}{2}\MM{u}\cdot M\MM{m}\,.
  $$
  Hence, proposition \ref{ham side} follows and we obtain the
  Hamiltonian (\ref{Ham-mom}) via the Legendre transform.
\end{proof}
Finally, we calculate Hamilton's canonical equations for this
Hamiltonian.
\begin{proposition}[Hamilton's canonical equations]
Hamilton's canonical equations with $H$ defined in equation 
(\ref{Ham-mom}) above may be expressed as
\begin{eqnarray*}
\MM{\dot{P}}_\beta & = & -\,\MM{P}_\beta
\cdot\sum_k\MM{u}_k\pp{\psi_k}{\MM{Q}_\beta}
(\MM{Q}_\beta)\,, \\
\MM{\dot{Q}}_\beta & = & \sum_k\MM{u}_k\psi_k(\MM{Q}_\beta)
\,.
\end{eqnarray*}
\end{proposition}
\begin{proof}
Hamilton's canonical equations are
$$
\MM{\dot{P}}_\beta = -\,\pp{H}{\MM{Q}_\beta}\,,
\qquad \MM{\dot{Q}}_\beta = \pp{H}{\MM{P}_\beta}\,.
$$
The $\MM{P}$ equation
projects onto grid variables as
\begin{eqnarray*}
\MM{\dot{P}}_\beta & = & -\,\pp{H}{\MM{Q}_\beta}
= -\sum_k\pp{H}{\MM{m}_k}\cdot \pp{\MM{m}_k}{\MM{Q}_\beta}
\\
& = & -\sum_{k}(\MM{P}_\beta\cdot\MM{u}_k)
\pp{\psi_k}{\MM{Q}_\beta}(\MM{Q}_\beta)\,. 
\end{eqnarray*}
Likewise, the $\MM{Q}$ equation projects onto grid variables via,
\begin{eqnarray*}
\MM{\dot{Q}}_\beta & = & \pp{H}{\MM{P}_\beta}
= \sum_k\pp{H}{\MM{m}_k}\cdot \pp{\MM{m}_k}{\MM{P}_\beta}
 \\
& = & \sum_{k}\MM{u}_k\psi_k(\MM{Q}_\beta)
\,,
\end{eqnarray*}
 as required. These are identical to equations (\ref{dotP}) and (\ref{dotQ}),
respectively.
\end{proof}

\subsection{Constructing a fully discrete method}
To construct a fully discrete method we use the standard variational
integrator approach as described in \cite{LeMaOrWe2003}, applied to
the constrained action principle in definition \ref{implicit action
principle}. We replace the integral over time by a Riemann sum over
discrete time levels, and define the map
\[
\phi:\Omega^{n_p}\times\Omega^{n_p}\to T\Omega^{n_p}
\]
which approximates $\MM{\dot{Q}}_\beta$. We write the discrete
action
\[
\mathcal{A}_d = \Delta t\sum_{n=0}^N
\left(L(\MM{u}^n) - \sum_\beta\MM{P}^n_\beta\cdot\phi(\MM{Q}^n,\MM{Q}^{n-1})
\right).
\]
Minimisation of the discrete action over $\MM{u}$, $\MM{P}$ and
$\MM{Q}$ gives the numerical scheme.

For example, consider the choice
\[
\phi(\MM{Q}^n,\MM{Q}^{n-1})=\frac{\MM{Q}^n-\MM{Q}^{n-1}}{\Delta t}.
\]
In this case, the discrete action becomes
\begin{eqnarray*}
\mathcal{A}_d & = & \Delta t\sum_{n=0}^N\left(
\sum_{kl}H_{kl}\MM{u}_k^n\cdot\MM{u}_l^n\right. \\
& & \quad + \left.
\sum_\beta\MM{P}_\beta^n\cdot\left(
\frac{\MM{Q}^n_\beta-\MM{Q}^{n-1}_\beta}{\Delta t}
-\sum_k\MM{u}_k^n\psi_k(\MM{Q}^{n-1}_\beta)\right)
\right),
\end{eqnarray*}
which is minimised by the solutions
\begin{eqnarray}
\label{1st order 1} \sum_lH_{kl}\MM{u}_l^n & = & 
\sum_\beta\MM{P}_\beta^n\psi_k(\MM{Q}_\beta^{n-1})\,, \\
\label{1st order 2} \MM{Q}^{n+1}_\beta & = & \MM{Q}^n_\beta
+\Delta t\sum_k\MM{u}_k^{n+1}\psi_k(\MM{Q}_\beta^n)\,, \\
\label{1st order 3} \MM{P}^{n+1}_\beta & = & \MM{P}^n_\beta
-\Delta t\MM{P}^{n+1}_\beta\cdot\sum_k\MM{u}^{n+1}
\pp{\psi_k}{\MM{Q}}(\MM{Q}^n_\beta)
\,.
\end{eqnarray}
This system is equivalent to the 1st order symplectic Euler-A method
(\emph{i.e.}  the 1st order symplectic method which is implicit in
$\MM{P}$ and explicit in $\MM{Q}$) applied to the Hamiltonian system
given in proposition \ref{ham side}.

\section{The discrete Euler-Poincar\'e equation for VPM}\label{DEP-sec}

In this section we compute the discrete EPDiff equation directly on
the Eulerian grid.
\rem{Need more intro: where are we going, where have we been, how does
our route so far prepare us to take the next step?}

\begin{theorem}[Discrete Euler-Poincar\'e theorem]
  With the above notation and assumptions, let $L(\MM{Q},\MM{\dot{Q}})$ be
  a Lagrangian expressible 
  as a function $L(\MM{u}(\MM{Q},\MM{\dot{Q}}))$ of grid velocity
representative $\MM{u}$ only.  The
  following four statements are equivalent:
  \begin{description}
    \item [(i)]
    The VPM trajectory $\MM{Q}(t)$ is an extremal of the constrained
    action
    \[
      {S} = \int_a^b L(\MM{u}) 
      + \sum_\beta\MM{P}_\beta\cdot
      \left(\MM{\dot{Q}}_\beta-\sum_k\MM{u}_k\psi_k(\MM{Q}_\beta)\right)
      \diff{t}
      \]
      with boundary conditions $\MM{Q}(a)=\MM{Q}_a$, $\MM{Q}(b)
      =\MM{Q}_b$, and where $\MM{P}_\beta$, $\beta=1,\ldots,n_p$, are
      Lagrange multipliers.
    \item [(ii)] The VPM trajectory $\MM{Q}(t)$ is the solution to the
      canonical Hamiltonian system in proposition \ref{ham side} with
      suitable boundary conditions.
    \item  [(iii)] The grid velocity $\MM{u}$ obtained from the grid
      momentum
      $\MM{m}$ using equation (\ref{discrete elliptic}), itself
      obtained from the particle momentum $\MM{P}$ using lemma
      (\ref{gridmom-eqn}), minimises the action ${S} =
      \int_a^b L(\MM{u})\diff{t}$
      upon taking constrained variations $\delta\MM{u}_k$ which satisfy
      $$
      [\delta\MM{u}]_\beta^P \equiv
      \sum_k\delta\MM{u}_k(t)\psi_k(\MM{Q}_\beta)
      = \sum_k\psi_k(\MM{Q}_\beta)\MM{\dot{w}}_k
      -[\ad_{\MM{u}}\MM{w}]_\beta^P,  
      =
      [\MM{\dot{w}}]_\beta^P -[\ad_{\MM{u}}\MM{w}]_\beta^P    
      \,,
      $$ 
      where $[\ad_{\MM{u}}\MM{w}]_\beta^P$ is defined by
      $$ 
      (\ad_{\MM{u}}\MM{w})_\beta \equiv 
      \sum_{k,\,l}\psi_k(\MM{Q}_\beta)
        \Big(\MM{u}_k(t)\cdot\pp{\psi_l}{\MM{Q}_\beta}
        (\MM{Q}_\beta)\Big)
        \MM{w}_l(t)-
        \psi_k(\MM{Q}_\beta)
	\Big(\MM{w}_k(t)\cdot\pp{\psi_l}{\MM{Q}_\beta}(\MM{Q}_\beta)
        \Big)\MM{u}_l(t)
      \,.      
      $$
      \item  [(iv)] The grid momentum $\MM{m} =
      (\MM{m}_1,\ldots,\MM{m}_{n_g})
    \in \mathbb{R}^{d\times n_g}$, satisfies the discrete
      Euler-Poincar\'e equation
      $$
      \MM{\dot{m}}_k + (\ad^*_{\MM{u}}\MM{m})_k = 0,
      \quad k=1,\dots,n_g\,,
      $$ 
      where
      \begin{eqnarray*}
        (\ad^*_{\MM{u}}\MM{m})_k &=& \sum_n(M^{-1})_{kn}
        \Big(\sum_l\sum_\beta\pp{\psi_l}{\MM{Q}_\beta}(\MM{Q}_\beta)
        (\MM{P}_\beta\cdot\MM{u}_l)
        \psi_n(\MM{Q}_\beta)
        \\&&
        \quad-\sum_l\delta_{nl}\sum_\beta\MM{P}_\beta
        \pp{\psi_l}{\MM{Q}_\beta}(\MM{Q}_\beta)
        \cdot\sum_m
\MM{u}_m(t)\psi_m(\MM{Q}_\beta)\Big),
      \end{eqnarray*}
      so that
      $$
      \langle \ad^*_{\MM{u}}\MM{m},\MM{w}\rangle_g = 
      \langle\MM{P},[\ad_{\MM{u}}\MM{w}]^P\rangle_p,
      $$
      where $\langle\cdot,\cdot\rangle_g$ is the grid inner product 
      defined by
      $$
      \langle \MM{f},\MM{g}\rangle_g 
      = \sum_{k,l}\MM{f}_k\cdot M_{kl}\MM{g}_l,
      $$
      (\emph{i.e.} a discrete approximation of the $L^2$ inner product in
      the 
      continuous case), $\langle\cdot,\cdot\rangle_p$ is the particle
      inner product on $T\Omega^{n_p}$ defined by
      $$
      \langle \MM{F},\MM{G}\rangle_p = \sum_{\beta}\MM{F}_\beta
      \cdot\MM{G}_\beta,
      $$
      and where $\MM{P}_\beta$ satisfies 
      \begin{equation}
\label{m,P EP}
      \sum_lM_{kl}\MM{m}_l = \sum_\beta\MM{P}_\beta\psi_k(\MM{Q}_\beta).
      \end{equation}
  \end{description}
\label{EP Theorem}
\end{theorem}
\begin{remark}
      The operation $(\ad_{\MM{u}}\MM{w})_\beta$ in {\bf(iii)} is the Lie
bracket among vector fields evaluated at the particle location
$\MM{Q}_\beta$.  The operation $\ad^*_{\MM{u}}\MM{m}$ is its dual with
respect to the pairing $\langle\cdot\,,\,\cdot\rangle_g$ on the grid.
\end{remark}
\begin{proof}
  {\bf(i)} $\Leftrightarrow$ {\bf(ii)} follows 
  from the proposition \ref{ham side}. \\

  To prove {\bf(i)} $\Leftrightarrow$ {\bf(iii)} we note that the
  constrained variational principle given in proposition \ref{implicit
  action principle} is equivalent to the Lagrange-d'Alembert principle
  \[
  \left(
  \dd{}{t}\pp{L}{\MM{\dot{Q}}_\beta}
  -\pp{L}{\MM{Q}_\beta}
  \right)\cdot\delta\MM{Q}_\beta=0,
  \]
  with constrained variations $\delta\MM{Q}\in D_{\VPM,\MM{Q}}$ and
  the constraint $\MM{\dot{Q}}\in D_{\VPM,\MM{Q}}$. The variations
  $\delta\MM{u}$ must be expressed in terms of the variations
  $\delta\MM{Q}$ and $\delta\MM{\dot{Q}}$ which follows by taking
  variations in equation (\ref{vector field def}):
  \begin{equation}
    \label{time then variations}
    \delta\MM{\dot{Q}}_\beta = \sum_k\delta\MM{u}_k(t)\psi_k(\MM{Q}_\beta)
    +\sum_k\MM{u}_k(t)\pp{\psi_k}{\MM{Q}_\beta}(\MM{Q}_\beta)\cdot\delta
    \MM{Q}_\beta.
  \end{equation}
  The variations $\delta\MM{Q}_\beta$ are written 
  \[
  \delta\MM{Q}_\beta = \sum_k\MM{w}_k(t)\psi_k(\MM{Q}_\beta),
  \]
  for some time series of velocity vectors on the grid $\MM{w}_k(t)$
  which vanishes on the end points.  Differentiating in time gives
  \begin{equation}
    \label{variations then time}
    \delta\MM{\dot{Q}}_\beta = \sum_k\MM{\dot{w}}_k(t)\psi_k(\MM{Q}_\beta)
    +\sum_k\MM{w}_k(t)\pp{\psi_k}{\MM{Q}_\beta}(\MM{Q}_\beta)\cdot
    \MM{\dot{Q}}_\beta.
  \end{equation}
  Combining equations (\ref{time then variations}) and (\ref{variations
  then time}) gives 
  \[
  \sum_k\delta\MM{u}_k(t)\psi_k(\MM{Q}_\beta)
  = \sum_k\psi_k(\MM{Q}_\beta)\left(
  \MM{\dot{w}}_k(t)+\MM{u}_k(t)\cdot
  \sum_l\pp{\psi_l}{\MM{Q}_\beta}(\MM{Q}_\beta)
  \MM{w}_l(t)-
  \MM{w}_k(t)\cdot\sum_l\pp{\psi_l}{\MM{Q}_\beta}(\MM{Q}_\beta)
  \MM{u}_l(t)
  \right),
  \]
  which we denote as
  \[
    [\delta\MM{u}]_\beta^P =
    \sum_k\delta\MM{u}_k(t)\psi_k(\MM{Q}_\beta)
    = \sum_k\psi_k(\MM{Q}_\beta)\MM{\dot{w}}_k-[\ad_{\MM{u}}\MM{w}]_\beta^P
    = [\MM{\dot{w}}]_\beta^P-[\ad_{\MM{u}}\MM{w}]_\beta^P
    \,.
    \]
This proves {\bf(i)}$\Leftrightarrow${\bf(iii)} and defines
$[\ad_{\MM{u}}\MM{w}]_\beta$.
\begin{remark}
  The bracket-subscript notation $[\,\cdot\,]_\beta^P$ introduced in
  the last formula emphasizes the VPM distinction between particle
  vector fields such as $[\MM{\dot{w}}]_\beta^P$ and their grid
  representatives $\MM{\dot{w}}_k$, related by
  $[\MM{\dot{w}}]_\beta^P=\sum_k\MM{\dot{w}}_k\psi_k(\MM{Q}_\beta)$.
  For example,
\begin{eqnarray*}
\langle \MM{m},\delta\MM{u}\rangle_g
\equiv
\sum_{k,l}\MM{m}_kM_{kl}\cdot\delta\MM{u}_l
=
\sum_{\beta,l}\psi_l(\MM{Q}_\beta)\MM{P}_\beta\cdot\delta\MM{u}_l
=
\sum_\beta\MM{P}_\beta\cdot[\delta\MM{u}]_\beta
\equiv
\langle \MM{P},[\delta\MM{u}]\rangle_p
\end{eqnarray*}
where the momentum relation (\ref{gridmom-eqn}) was used in the second
step and the relation between VPM particle vector fields and their
grid representatives was applied in the third step. A similar
calculation allows one to write the dual relations defining the VPM
particle- and grid-representatives of ad$^*$. Namely,
\begin{eqnarray*}
\langle \MM{P},[\ad_{\MM{u}}\MM{w}]\rangle_p
=
\langle \MM{m},\ad_{\MM{u}}\MM{w}\rangle_g
\,,
\end{eqnarray*}
whose dual relation may be conveniently written as
\begin{eqnarray*}
\langle [\ad^*_{\MM{u}}\MM{P}]^G,\MM{w}\rangle_p
=
\langle \ad^*_{\MM{u}}\MM{m},\MM{w}\rangle_g
\,,
\end{eqnarray*}
in order to define ad$^*$ in both particle and grid representations. 
In particular, this implies 
      $$
      \langle \ad^*_{\MM{u}}\MM{m},\MM{w}\rangle_g = 
      \langle\MM{P},[\ad_{\MM{u}}\MM{w}]^P\rangle_p\,,
      $$
as claimed in the theorem. 
\end{remark}
To prove {\bf(iii)}$\Leftrightarrow${\bf(iv)} we take variations 
$\delta\MM{u}$ in ${S}$:
\begin{eqnarray*}
0=\delta{S}
& = & \int
\langle \MM{m},\delta\MM{u}\rangle_g
\,dt
\\
& = & \int
\langle \MM{P}, [\delta\MM{u}] \rangle_\beta
\,dt
\\
& = & \int
\langle\MM{m},
\MM{\dot{w}}\rangle_g-\langle\MM{P},[\ad_{\MM{u}}\MM{w}]\rangle_\beta
\,dt
\\
& = & \int
\langle-\MM{\dot{m}},\MM{w}\rangle_g
-\sum_\beta\MM{P}_\beta
\cdot\left(
\MM{u}_k(t)\pp{\psi_k}{\MM{Q}_\beta}(\MM{Q}_\beta)\cdot
\sum_l\MM{w}_l(t)\psi_l(\MM{Q}_\beta)\right. \\
& & \hspace{2cm}
-\left.\MM{w}_k(t)\pp{\psi_k}{\MM{Q}_\beta}(\MM{Q}_\beta)\cdot
\sum_l\MM{u}_l(t)\psi_l(\MM{Q}_\beta)\right)
\,dt \\
& = & \int
\langle-\MM{\dot{m}},\MM{w}\rangle_g
-\sum_k
\MM{w}_k(t)\cdot\left(
\sum_{l,\beta}\pp{\psi_l}{\MM{Q}_\beta}(\MM{Q}_\beta)
(\MM{P}_\beta\cdot\MM{u}_l)
\psi_k(\MM{Q}_\beta)
\right. \\
& & \hspace{2cm}
\left.
-
\sum_{l,\beta}\delta_{kl}\MM{P}_\beta
\pp{\psi_l}{\MM{Q}_\beta}(\MM{Q}_\beta)
\cdot\sum_m
\MM{u}_m(t)\psi_m(\MM{Q}_\beta)
\right)
\,dt \\
& = & -\int
\langle\frac{d}{dt}\MM{m}
+\ad^*_{\MM{u}}\MM{m}
,\MM{w}\rangle_g
\,dt
\end{eqnarray*}
where we have integrated by parts. The grid representation $\MM{w}$ is
arbitrary and therefore
$$
\frac{d}{dt}\MM{m}
+\ad^*_{\MM{u}}\MM{m}=0,
$$ as required. 
\end{proof}
The correspondence {\bf(ii)}$\Leftrightarrow${\bf(iv)} was also proved
by direct calculation in \cite{VPM}.


\section{Left action momentum map}\label{L-Action-sec}

First, recall that a canonical action of a Lie algebra $A$ on a
symplectic manifold $\mathcal{M}$ is a mapping from $A$ to Hamiltonian
vector fields on $\mathcal{M}$ which preserves the Lie
brackets. Consider an element $\MM{\xi}$ of $A$ and its action
$\MM{\xi}_{\mathcal{M}}$ on $\mathcal{M}$ which has Hamiltonian $J$.
The momentum map $\MM{J}$ is related to the Hamiltonian $J$ by
$$
\langle\MM{J},\MM{\xi}\rangle = J,
$$
for all such elements $\MM{\xi}$, where
$\langle\,\cdot\,,\,\cdot\,\rangle:A^*\times A\to\mathbb{R}$ is the inner
product between $A$ and and its dual $A^*$.

If $A$ acts on a manifold $\mathcal{M}$ then we can define a canonical
action of $A$ on $T^*\mathcal{M}$ with Hamiltonian
$$
J = \langle(\MM{P},\MM{Q}),\MM{\xi}_{\mathcal{M}}\rangle
= \MM{P}\cdot\MM{\xi}_{\mathcal{M}}.
$$
This is called the {\bfi cotangent lift} of the action to
$T^*\mathcal{M}$. The definition of the momentum map for the cotangent
lift of an action then becomes
\begin{equation}
\label{cotangent mom map}
\langle\MM{J},\MM{\xi}\rangle = 
\langle(\MM{P},\MM{Q}),\MM{\xi}_{\mathcal{M}}\rangle.
\end{equation}

We define the left-action of $\mathfrak{X}(\Omega)$ on $\Omega$ 
by
\[
\MM{\xi} \mapsto \MM{\xi}_{\MM{Q}}
=\MM{\xi}(\MM{Q})\cdot\pp{}{\MM{Q}}.
\]
The Hamiltonian for the cotangent-lifted left action is then
$$
J(\MM{\xi})(\MM{P,Q}) = \langle(\MM{P},\MM{Q}),\MM{\xi}_{\MM{Q}}\rangle.
$$ We wish to obtain a momentum map which maps into the representation
of $D_{\VPM}$ given by the map $\mathbb{R}^{d\times
n_g}\to\mathfrak{X}(\Omega)$:
\[
\MM{u} \mapsto \sum_k\MM{u}_k\psi_k(\MM{Q}).
\]
We do this by restricting $J(\MM{\xi})$
to elements of $D_{\VPM}$:
\begin{eqnarray*}
J(\MM{u})(\MM{P},\MM{Q})
&=& \left\langle(\MM{P},\MM{Q}),\sum_k\MM{u}_k\psi_k(\MM{Q})\cdot\pp{}{\MM{Q}}
\right\rangle, \\
& = & \sum_{\beta k}\MM{P}_\beta\cdot\MM{u}_k\psi_k(\MM{Q}),
\end{eqnarray*}
and this relation defines the left action momentum map
\begin{equation}\label{mommap}
\MM{J}_k^L(\MM{P},\MM{Q})=\sum_\beta\MM{P}_\beta\psi_k(\MM{Q}_\beta)\,.
\end{equation}
As mentioned in Remark \ref{left mom-map} this is again equation
(\ref{gridmom-eqn}) derived earlier from constrained variations of the
VPM action (\ref{VPM action}) with respect to the grid representatives of
the velocity. This momentum map is the discrete version for VPM of a
general result for Clebsch variational principles for ideal fluid
dynamics \cite{HoKu1983,HoMaRa1998}. 


\section{Right action momentum map}\label{R-Action-sec}

\paragraph{Right action of $\mathfrak{X}(\Omega)$ on
$\{\MM{Q}_\beta(t)\}_{\beta=1}^{n_p}$}

Next we define the right action of $\Diff(\Omega)$. To do this we
require the entire ``history'' of $\MM{u}(t)$. Given initial
conditions $\MM{Q}(0)=\MM{Q}^0$, the history of $\MM{u}(t)$ produces a
solution
$\MM{Q}(t)$ with
$$
\MM{\dot{Q}}(t) = \sum_k\MM{u}_k(t)\psi_k(\MM{Q}(t))\,, \qquad
\MM{Q}(0) = \MM{Q}^0.
$$
This solution can be extended to a one-parameter family of
diffeomorphisms $g_t:\Omega\to\Omega$ with
\begin{equation}
\pp{}{t}g_t(x) = \sum_k\MM{u}_k(t)\psi_k(g_t(x))\,, \qquad
g_0(x) = x\,. \label{g def}
\end{equation}
In particular, $g_t(\MM{Q}^0_\beta) = \MM{Q}_\beta(t)$. This allows
a right action of $\eta\in\Diff(\Omega)$ on $g_t$ to be defined
\emph{via} composition:
$$
g_t\mapsto g_t\cdot\eta = g_t\circ\eta\,.
$$
Using the tangent map, one may define a right action of
$\MM{\xi}\in\mathfrak{X}(\Omega)$ on $g_t$:
$$
g_t\mapsto g_t\cdot\MM{\xi} = \pp{g_t}{x}\cdot\MM{\xi}\,.
$$
Again in particular,
$$
\MM{Q}_\beta(t) \mapsto \pp{g_t}{x}(\MM{Q}^0)\cdot\MM{\xi}(\MM{Q}^0_\beta)
\,.
$$
Differentiating equation (\ref{g def}) gives
\begin{equation}
\label{g evol}
\pp{}{t}\pp{g_t}{x}(\MM{Q}^0_\beta) = 
\sum_k\MM{u}_k(t)\pp{\psi_k(g_t(\MM{Q}^0_\beta))}
{x}\cdot\pp{g_t}{x}(\MM{Q}^0_\beta) = \sum_k\MM{u}_k(t)\pp{\psi_k}
{x}(\MM{Q}_\beta(t))\cdot\pp{g_t}{x}(\MM{Q}^0_\beta)\,,
\end{equation}
with initial conditions
$$
\pp{g_0}{x}(\MM{Q}^0_\beta) = \Id.
$$ 
This means that, given $\MM{u}_k(t)$ (and hence, given
$\{\MM{Q}_\beta(t)\} _{\beta=1}^{n_p}$), the Jacobian
${\partial g_t/\partial x}(\MM{Q}^0_\beta)$ may be obtained without needing
to calculate $g_t$ as a map over the whole of $\Omega$.

We can interpret this map as a canonical momentum map by extending the
canonical coordinates $\{\MM{P}_{\beta},\MM{Q}_{\beta}\}_{n=1}^{n_p}$
to $\{\MM{P}_{\beta},\MM{Q}_\beta,J_{\beta}\}_{n=1}^{n_p}$ where
$J_\beta$ is a $d\times d$ matrix for each $\beta$. We write $\MM{J}$
as a column vector, \emph{e.g.,} in two dimensions
\[
\MM{J} = \left(J_{11,1},J_{12,1},J_{21,1},J_{22,1},\ldots,
J_{11,n_p},J_{12,n_p},J_{21,n_p},J_{22,n_p}\right)^T,
\]
and consider
the Hamiltonian system, cf. \cite{HoKuLe1985}
\[
\begin{pmatrix}
\MM{\dot{Q}} \\
\MM{\dot{P}} \\
\MM{\dot{J}} \\
\end{pmatrix}
=
\begin{pmatrix}
0 & \Id & 0 \\
-\Id & 0 & -B^T(\MM{Q})K(\MM{Q},\MM{J})^T \\
0 & K(\MM{Q},\MM{J})B(\MM{Q}) & 0 \\
\end{pmatrix}
\begin{pmatrix}
\nabla_{\MM{Q}} \\
\nabla_{\MM{P}} \\
\nabla_{\MM{J}} \\
\end{pmatrix}
H\,,
\]
where
\[
B\MM{\dot{Q}} = \MM{u}, \quad\mbox{if}\quad \MM{\dot{Q}}_\beta
=\sum_k\MM{u}_k\psi_k(\MM{Q}_\beta)\,,
\]
and 
\[
\left(K(\MM{Q},\MM{J})\MM{u}\right)_{ij,\beta} 
= \sum_{kl}u_{k,i}\pp{\psi_k}{Q_l}(\MM{Q}_\beta)
J_{lj,\beta}\,.
\]
When the Hamiltonian is a function of $\MM{P}$ and $\MM{Q}$ only, we
recover the canonical Hamiltonian structure for $\MM{P}$ and $\MM{Q}$.
Furthermore, if the Hamiltonian is a function of grid momentum only,
so that
\[
\MM{\dot{Q}}_\beta = \sum_k\MM{u}_k\psi_k(\MM{Q}_\beta),
\]
for some $\MM{u}$, then
\[
\dot{J}_\beta = \sum_k\MM{u}_k\pp{\psi_k}{\MM{Q}}(\MM{Q}_\beta)
\cdot J_\beta,
\]
as required. This larger system enables us to talk about
discrete particle-relabelling, as summarised in the following
theorem:
\begin{theorem}\label{Inv-thm}
Consider the time-continuous VPM discretisation of EPDiff, given
in the above enlarged space, with Hamiltonian
\[
\mathcal{H} = \frac{1}{2}\sum_{klij}
M_{ik}\MM{m}_k\cdot(H^{-1})_{kl}M_{lj}\MM{m}_j,
\]
with $H$ the Helmholtz operator and $M$ the mass matrix. Then the 
flows of the vector field with Hamiltonian
\[
h = \sum_\beta\MM{P}_\beta\cdot J_\beta\cdot\MM{\xi}_\beta,
\]
for any constant vector $\{\MM{\xi}_\beta\}_{\beta=1}^{n_p}$,
leave the Hamiltonian $\mathcal{H}$ invariant.
\end{theorem}
\begin{proof}
The Hamiltonian $h$ generates the flow
\begin{eqnarray*}
\MM{\dot{Q}} & = & J_\beta\cdot\MM{\xi}, \\
\MM{\dot{P}} & = & -B^T(\MM{Q})K^T(\MM{Q},\MM{J})(\MM{P}\MM{\xi}),
\end{eqnarray*}
for $\MM{P}$ and $\MM{Q}$, where we write
\[
\left(\MM{J}\cdot\MM{\xi}\right)_\beta =
J_\beta\cdot\MM{\xi}_\beta, \quad \beta=1,\ldots,n_p,
\]
and
\[
\left(\MM{P}\MM{\xi}\right)_{ij,\beta} =
P_{i,\beta}\xi_{j,\beta}, \quad \beta=1,\ldots,n_p, \qquad i,j=1,\ldots,
d.
\]
The Hamiltonian $\mathcal{H}$ is a function of grid momentum $\MM{m}$ only
so it suffices to check that the quantity
\[
\sum_\beta\MM{P}_\beta\psi_k(\MM{Q}_\beta), \quad k=1,\ldots,n_g,
\] 
is invariant under the flow. We can check this directly, as
\begin{eqnarray*}
\dede{}{t}\sum_\beta\MM{P}_\beta\psi_k(\MM{Q}_\beta) & = &
\sum_\beta\left(\MM{\dot{P}}_\beta\psi_k(\MM{Q}_\beta) +
\MM{P}_\beta\pp{\psi_k}{\MM{Q}}(\MM{Q}_\beta)
\cdot\MM{\dot{Q}}_\beta\right), \\ & =
& \sum_\beta\left(
-\left(B^T(\MM{Q})A^T(\MM{Q},\MM{K})\MM{P}\MM{\xi}\right)_\beta
\psi_k(\MM{Q}_\beta) + 
\MM{P}_\beta\pp{\psi_k}{\MM{Q}}(\MM{Q}_\beta)\cdot
J_\beta\cdot\MM{\xi}_\beta
\right),
\end{eqnarray*}
and the first term becomes
\begin{eqnarray*}
-\sum_\beta\psi_k(\MM{Q}_\beta)
\left(B^T(\MM{Q})K^T(\MM{Q},\MM{J})\MM{P}\MM{\xi}\right)_\beta
& = & -\sum_\beta K^T_{k\beta}(\MM{Q},\MM{J})(\MM{P}\MM{\xi})_\beta,  \\
& = & -\sum_\beta\MM{P}_\beta\pp{\psi_k}{\MM{Q}}(\MM{Q}_\beta)
\cdot J_\beta\cdot\MM{\xi}_\beta,
\end{eqnarray*}
showing that the momentum is invariant, as required.
\end{proof}
\begin{corollary}\label{Right-mom-map}
The momentum map $\MM{J}^R$ with
\[
\MM{J}^R_\beta = \MM{P}_\beta\cdot J_\beta,
\]
is conserved for solution of semi-discrete EPDiff.
\end{corollary}
\begin{proof}
The result follows directly from Noether's theorem, {\it i.e.,} from
invariance of the Hamiltonian $\mathcal{H}$ in theorem \ref{Inv-thm}.
\end{proof}
\begin{remark}
The symmetry which changes $\MM{P}$ and $\MM{Q}$ while leaving
$\MM{m}$ invariant on the grid is our discrete form of the
particle-relabelling symmetry.  Next, we shall see that this symmetry
results in a discrete version of Kelvin's circulation theorem.
\end{remark}

\rem{Thus we have a dual-pair construction of momentum maps
corresponding to the left- and right actions, summarised in the
diagram \cite{MaWe1983}:

\begin{picture}(150,100)(-50,0)%
\put(120,75){$T^{\ast}\Omega^{n_p}$} 

\put(78,50){$\mathbf{J}^L$}        

\put(160,50){$\mathbf{J}^R$}   

\put(78,15){$D_{\VPM}^*$}

\put(160,15){$\mathfrak{X}(\Omega^{n_p})^{\ast}$}       

\put(130,70){\vector(-1, -1){40}}  

\put(135,70){\vector(1,-1){40}}  
\end{picture}}

\section{Kelvin's circulation theorem for discrete EPDiff}
\label{kelvin}
As discussed earlier, the discrete EPDiff Lagrangian is invariant
under the right action of $\mathfrak{X}(\Omega)$. This means that
$J^R$ in Corollary \ref{Right-mom-map} is a conserved momentum. In
particular,
$$
\dede{}{t}\Big(\MM{P}_\beta\cdot\pp{g_t}{x}(\MM{Q}^0_\beta)\Big)=0, 
\qquad \mbox{ no sum on }\beta
$$
for each $\beta=1,\ldots,n_p$. (This is obtained by integrating
$\MM{J}^R$ against a suitable function whose support contains only
$\MM{Q}_\beta$.)

We can interpret this result to prove a discrete form of Kelvin's
circulation theorem. Consider a loop $C(t)$ in $\Omega$ which is
embedded in the flow, \emph{i.e.},
$$
C(t) = g_t(C(0)).
$$
We choose $C(t)$ so that some of the particles with
trajectories $\MM{Q}_\beta(t)$ are located at the initial time $t=0$ on
$\MM{Q}_\beta(0)\in C(0)$. As
$g_t(\MM{Q}^0_\beta)=\MM{Q}_\beta(t)$, those particles will stay on
$C(t)$ for all time. Define the set $\Gamma$ so that $\beta\in\Gamma$
if $\MM{Q}_\beta^0$ is located on $C(0)$.

In order to discuss the circulation theorem, we need to introduce
a discretisation of density. As discussed in \cite{VPM}, this is 
done by associating a constant $\tilde{D}_\beta$ with each particle,
so that the density on the grid may be written 
\[
D_k = \sum_\beta\tilde{D}_\beta\psi_k(\MM{Q}_\beta).
\]
This allows us to represent $\MM{m}/D$ evaluated at the location of 
particle $\beta$ as $\MM{P}_\beta/\tilde{D}_\beta$. 

Next we need to approximate line integration round $C(t)$. We do this
by writing
\[
\int_{C(t)}\frac{\MM{m}}{D}\cdot\dt{\MM{x}} =
\int_0^{2\pi}\frac{\MM{m}\circ\MM{\gamma}_t}
{D\circ\MM{\gamma}_t}
\cdot \dede{\MM{\gamma_t}}{s}\dt{s},
\]
where $\gamma_t:[0,2\pi]\to C(0)$ is a parameterisation of the
loop $C(t)$. Substituting $\gamma_t = g_t\circ\gamma_0$ yields
\[
\int_{C(t)}\frac{\MM{m}}{D}\cdot\dt{\MM{x}} =
\int_0^{2\pi}\frac{\MM{m}\circ\MM{\gamma}_t}
{D\circ\MM{\gamma}_t}
\cdot\cdot\pp{g_t}{x}(\MM{\gamma_0}(s))
\cdot\dede{\MM{\gamma_0}}{s}\dt{s},
\]
which we can approximate with a Riemann sum
\[
\int_{C(t)}\frac{\MM{m}}{D}\cdot\dt{\MM{x}} \approx
\sum_\beta\frac{\MM{P}_\beta}{\tilde{D}_\beta}\cdot
\Delta x_\beta,
\]
where
\[
\Delta x_\beta =
J_\beta\cdot \dede{\MM{\gamma}_0}{s}
(s_\beta)\Delta s_\beta,
\]
with $\gamma_0(s_\beta)=\MM{Q}_\beta^0$, and $\Delta s_\beta
= s_{\beta+1}-s_{\beta}$. 

Using this discretised line integration scheme, we can state
our Kelvin circulation theorem as follows:
\begin{proposition}
Let $\{u_k(t)\}_{k=1}^{n_g}$ satisfy the discrete Euler-Poincar\'e
equations above, with $\{D_k(t)\}_{k=1}^{n_g}$ satisfying the discrete
density equation. Let $C(t)$ be a closed loop advected in the flow
generated by the velocity
$$\MM{u}(x,t)=\sum_k\MM{u}_k(t)\psi_\beta(\MM{Q}_\beta),$$ 
containing some subset of particles $Q_{\beta}$,
with $\beta\in B\subset(1,\ldots,n_p)$.  Define the discrete
circulation sum \\
$$
I(t) = \sum_{\beta\in B}\frac{\MM{P}_{\beta}}
{\tilde{D}_\beta}\cdot \Delta\MM{x}_\beta,
$$
Then $I(t)$ satisfies
\[
\frac{\diff}{\dt{t}}
I(t) = 0.
\]
\end{proposition}
\begin{proof}
The result proceeds directly from corollary \ref{Right-mom-map} for the
right action momentum map, which satisfies 
\[
\frac{d}{dt}(\MM{P}_\beta\cdot J_\beta)=0,
 \quad \forall \beta.
\]
\end{proof}

\section{Numerical Results}
\label{numerics}
\begin{figure}[htp]
\begin{center}
\scalebox{0.35}{\includegraphics{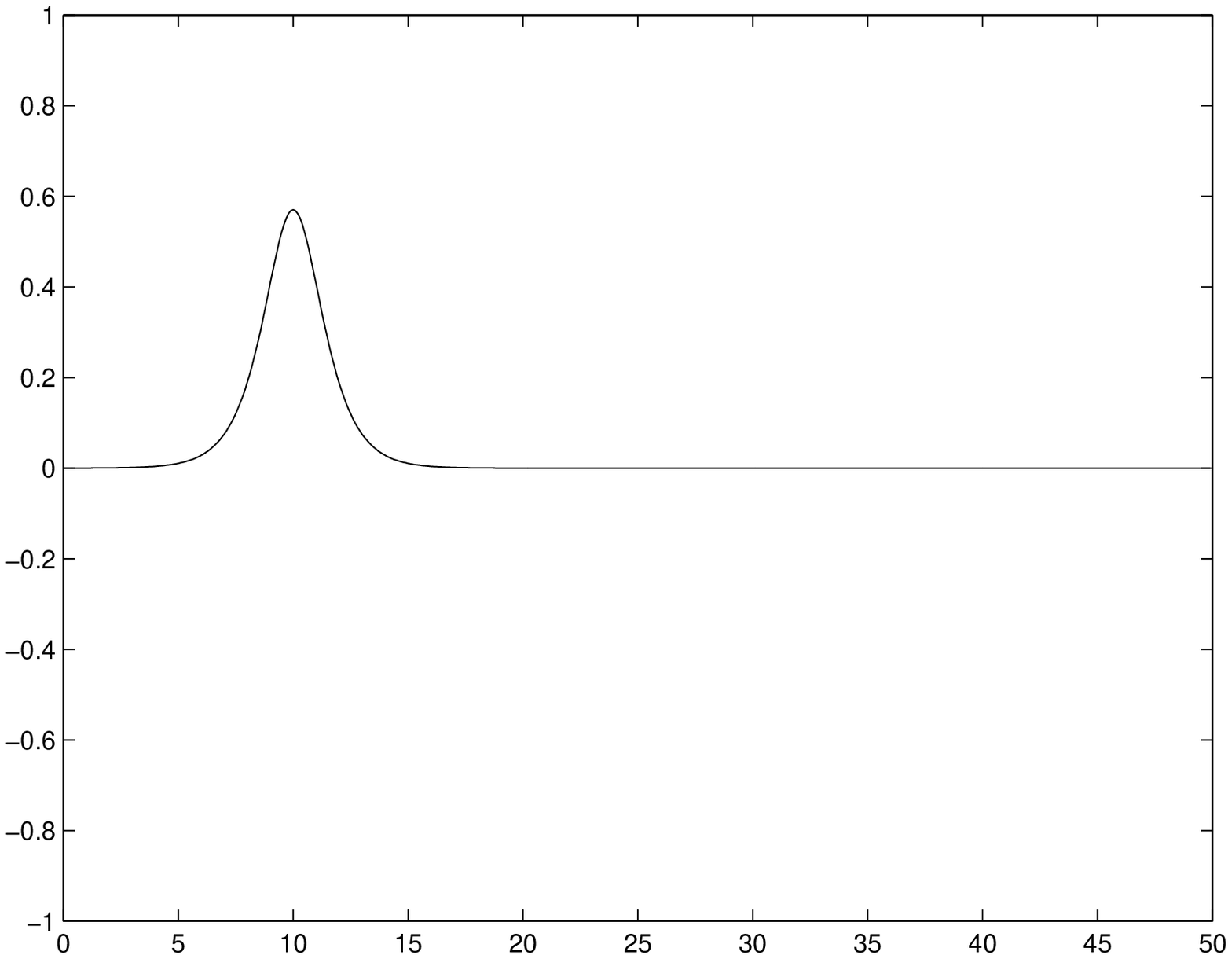}}
\scalebox{0.35}{\includegraphics{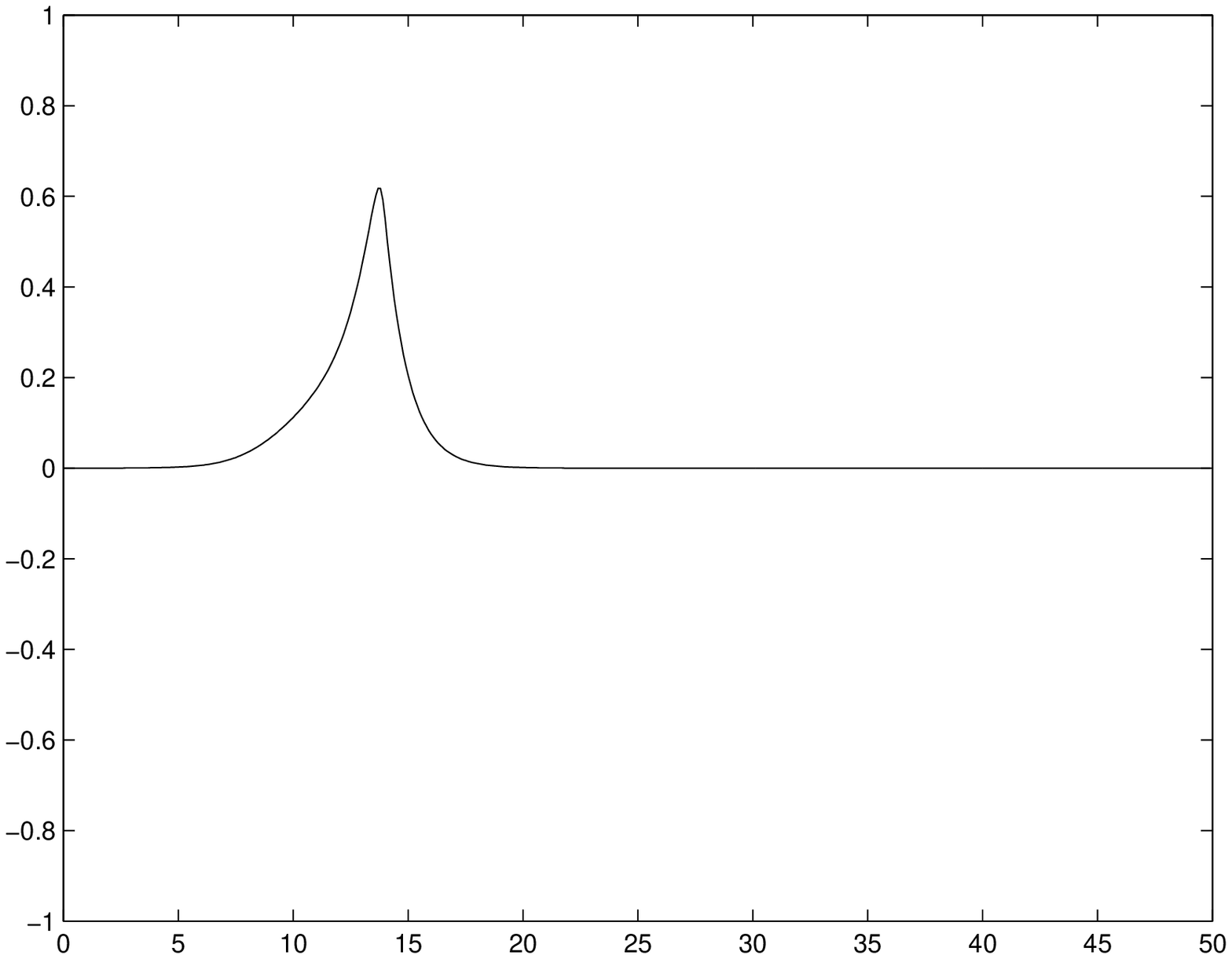}} \\
\scalebox{0.35}{\includegraphics{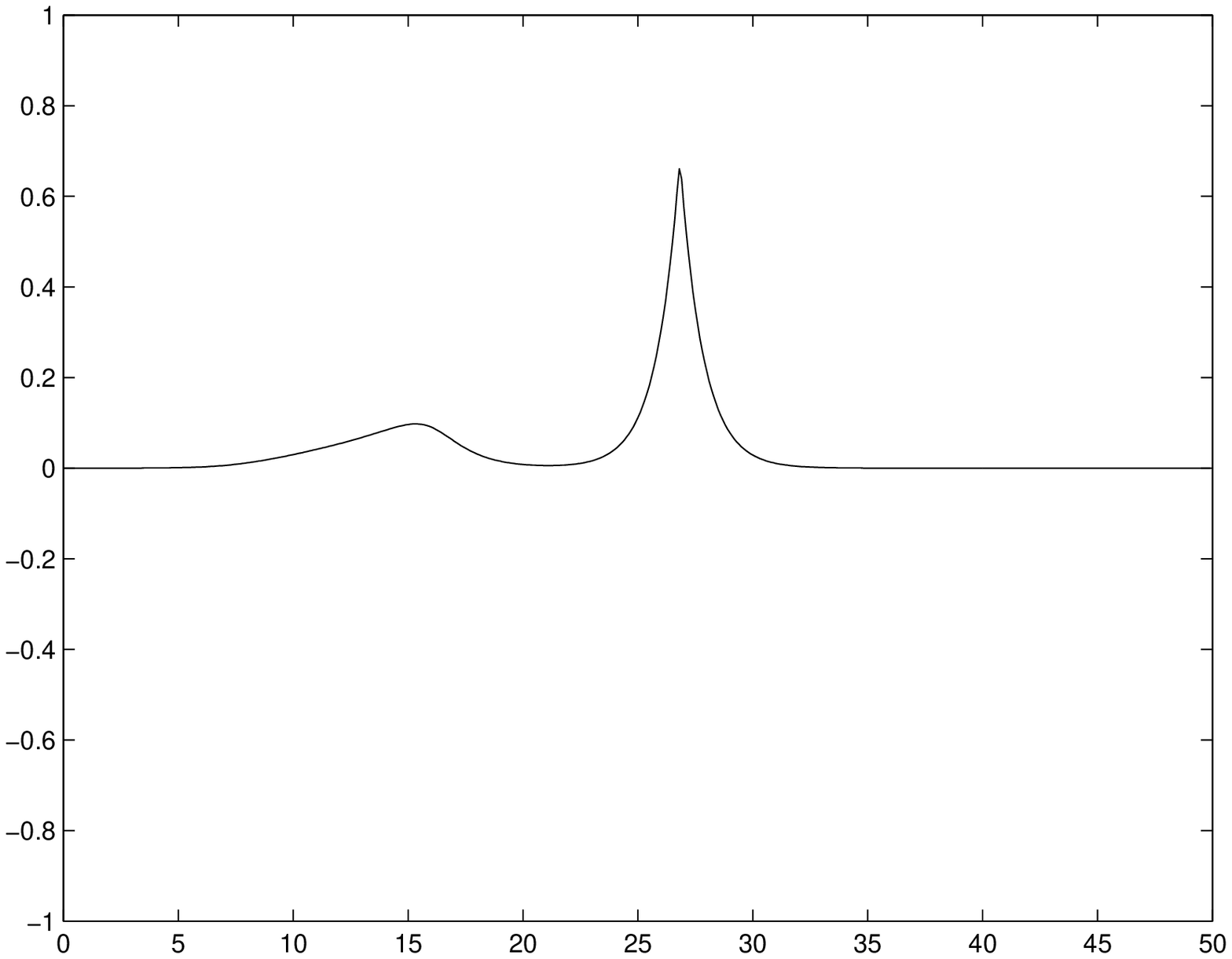}}
\scalebox{0.35}{\includegraphics{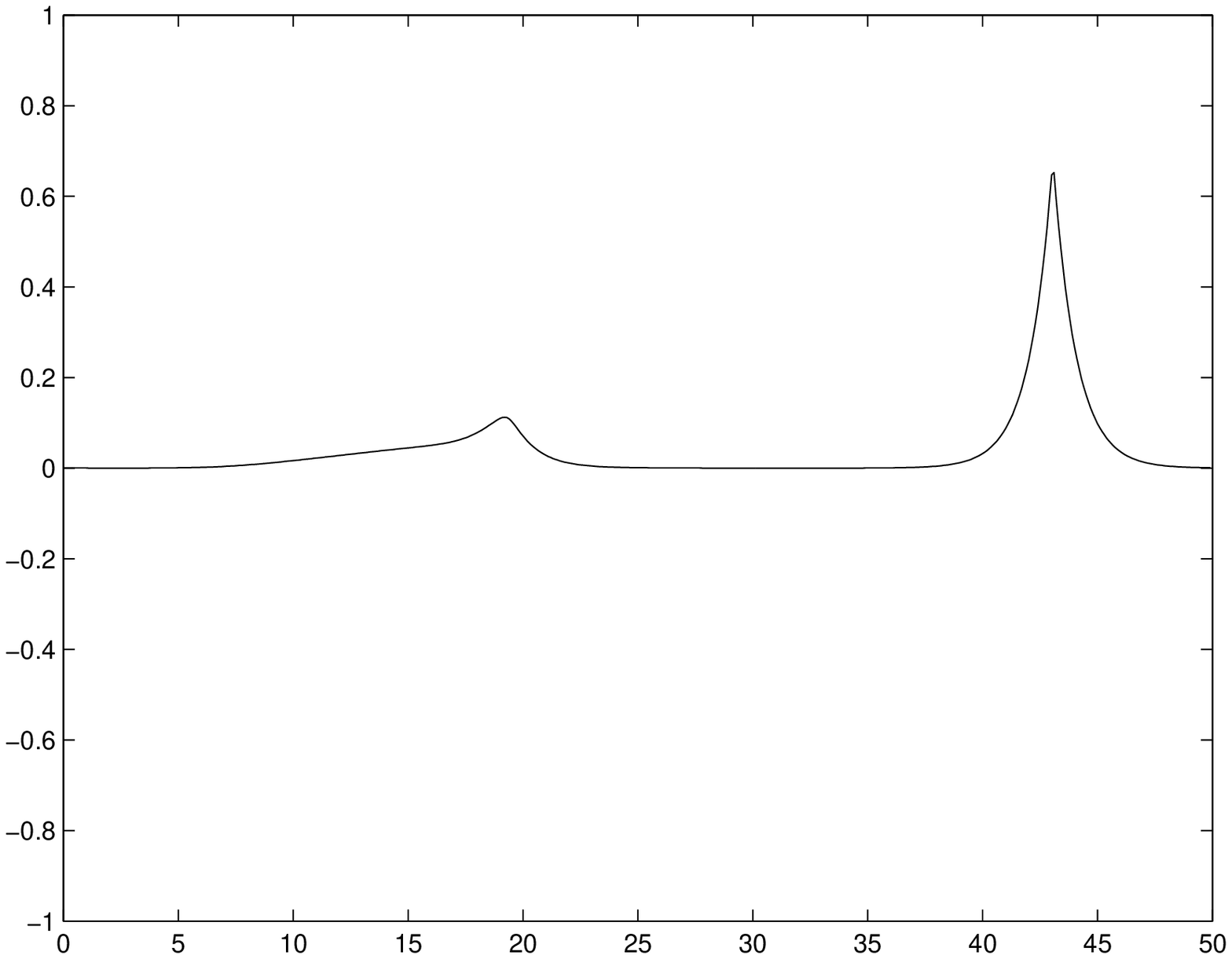}} \\
\end{center}
\caption{\label{sech_initial_value}Numerical solution to the 1D
  EPDiff equation with initial data $u=(\pi/2)e^x-2\sinh
  x\arctan(e^x)-1$, and scaling constant $\alpha=1$. The solutions are
  obtained using the VPM method with 500 grid points, 1000 particles,
  cubic B-splines for the basis functions, linear finite elements for
  the grid discretisation of the Lagrangian and a timestep of 0.1 with
  the first order time discretisation given by equations (\ref{1st order 1}
  -\ref{1st order 3}). The figures show the velocity field
  at times 0, 5, 25 and 50. At $t=5$ the smooth initial condition has
  ``leaned to the right'' and a discontinuous peak has formed. By
  $t=25$ the peakon is well separated from the smooth part of the
  solution and by $t=50$ a second peak is starting to form.}
\end{figure}

{\bfseries Convergence tests}
We begin by performing a convergence test for the 1D equations
\[
m_t + um_x + 2mu_x = 0, \quad (1-\alpha^2\partial_x^2)u,
\]
which, as discovered in \cite{CaHo1993}, is completely integrable with
the initial value problem dominated by peaked solitons (peakons) whose
first derivatives are discontinuous. This property is illustrated in
figure \ref{sech_initial_value} which shows a numerical integration 
of the 1D equations starting from a smooth initial condition with
singular peaked solitons emerging in finite time. 

\begin{figure}[htp]
\begin{center}
\scalebox{0.5}{\includegraphics{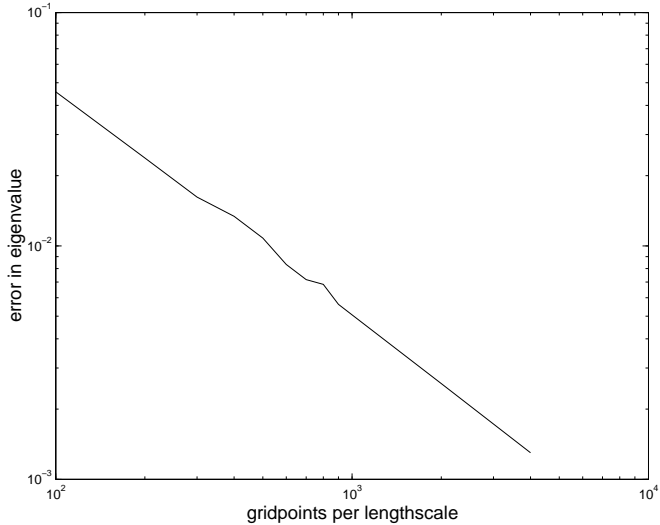}}
\end{center}
\caption{\label{eigplot}
  Plot of error in calculating the asymptotic
  speed of the first emitted peakon from an initial condition
  $u=(\pi/2)e^x-2\sinh x\arctan(e^x)-1$ against grid resolution
  (measured as number of gridpoints in one characteristic length
  $\alpha=1$), using the same method as figure
  \ref{sech_initial_value}.  The number of particles used was twice
  the number of grid points, and the timestep used was scaled with the
  grid size to guarantee convergence of the fixed-point method for
  solving the linear system (\emph{i.e.} not for accuracy). The
  logarithmic plot has slope -1 giving a linear scaling for the error
  with grid resolution.}
\end{figure}

For our first convergence test we use the result given in
\cite{CaHo1993} that for an initial condition $u=(\pi/2)e^x-2\sinh
x\arctan(e^x)-1$, with scaling constant $\alpha=1$, the asymptotic
speeds of the emitted peakons are $2/[(2n+1)(2n+3)]$,
$n=0,1,2,\ldots$. In particular the asymptotic speed of the first
peakon is $2/3$. Figure \ref{eigplot} shows that the numerical
calculation of the speed converges to the correct answer with a linear
scaling for error against grid resolution.

\begin{figure}[htp]
\begin{center}
\scalebox{0.45}{\includegraphics{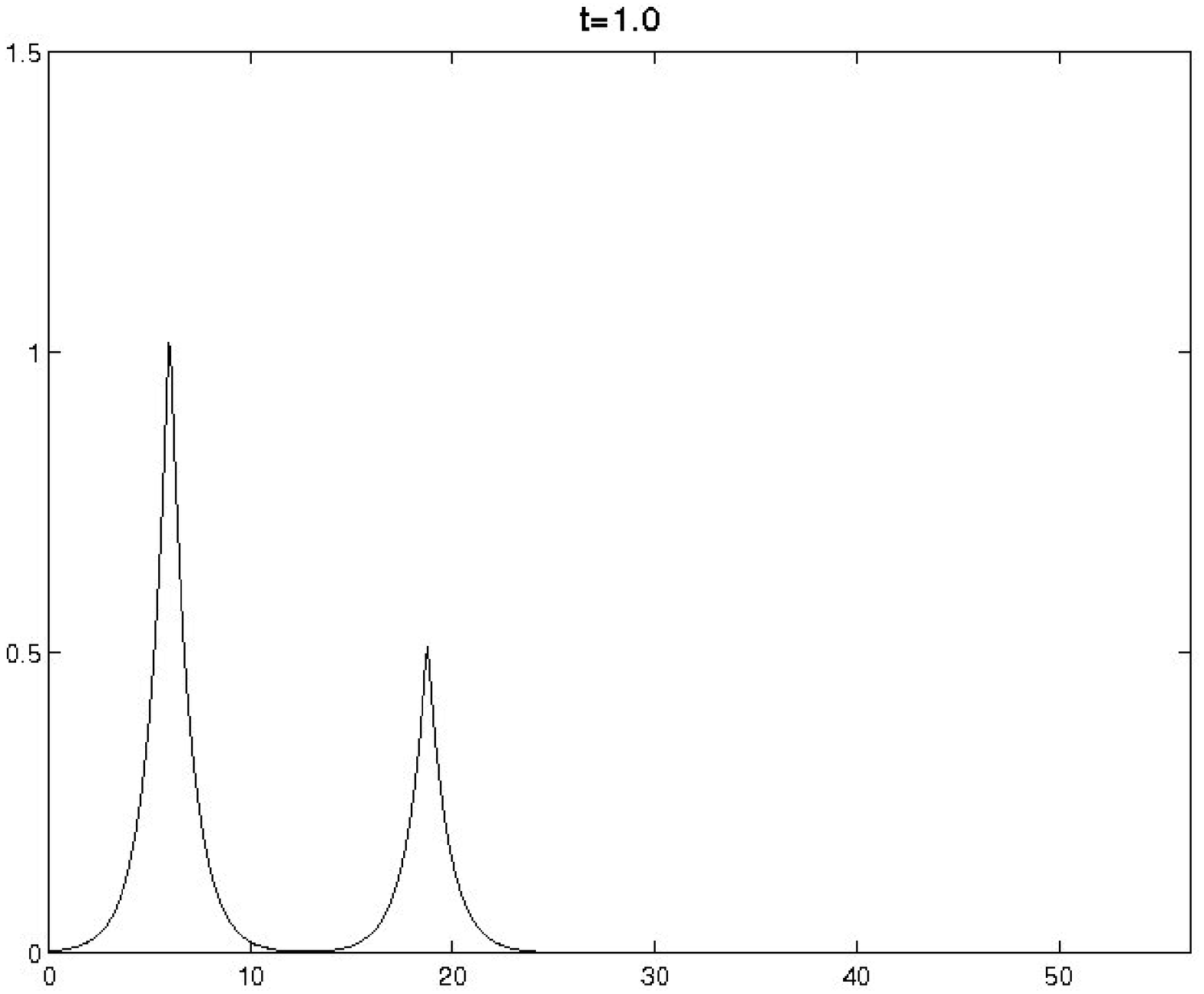}}
\scalebox{0.45}{\includegraphics{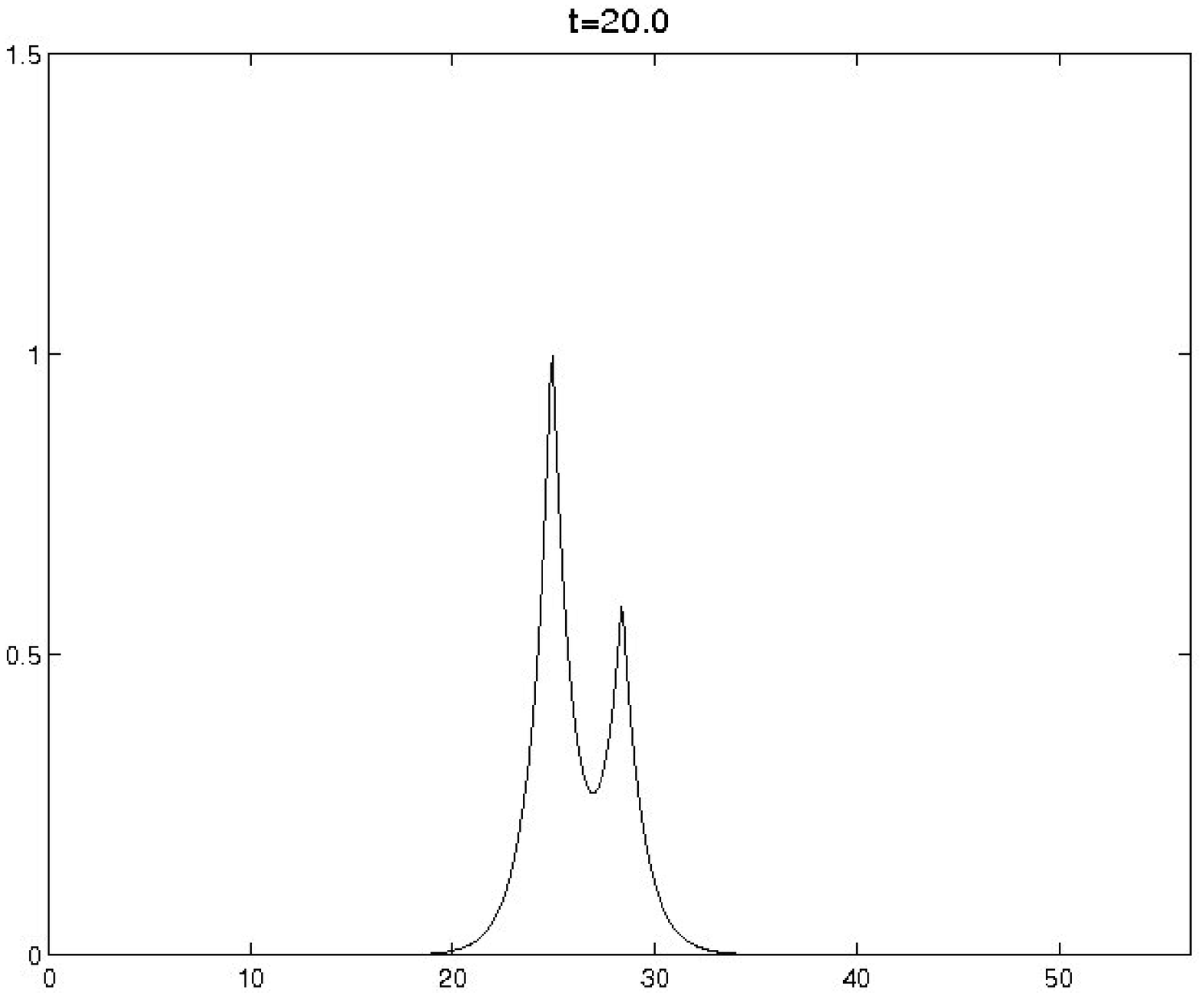}} \\
\scalebox{0.45}{\includegraphics{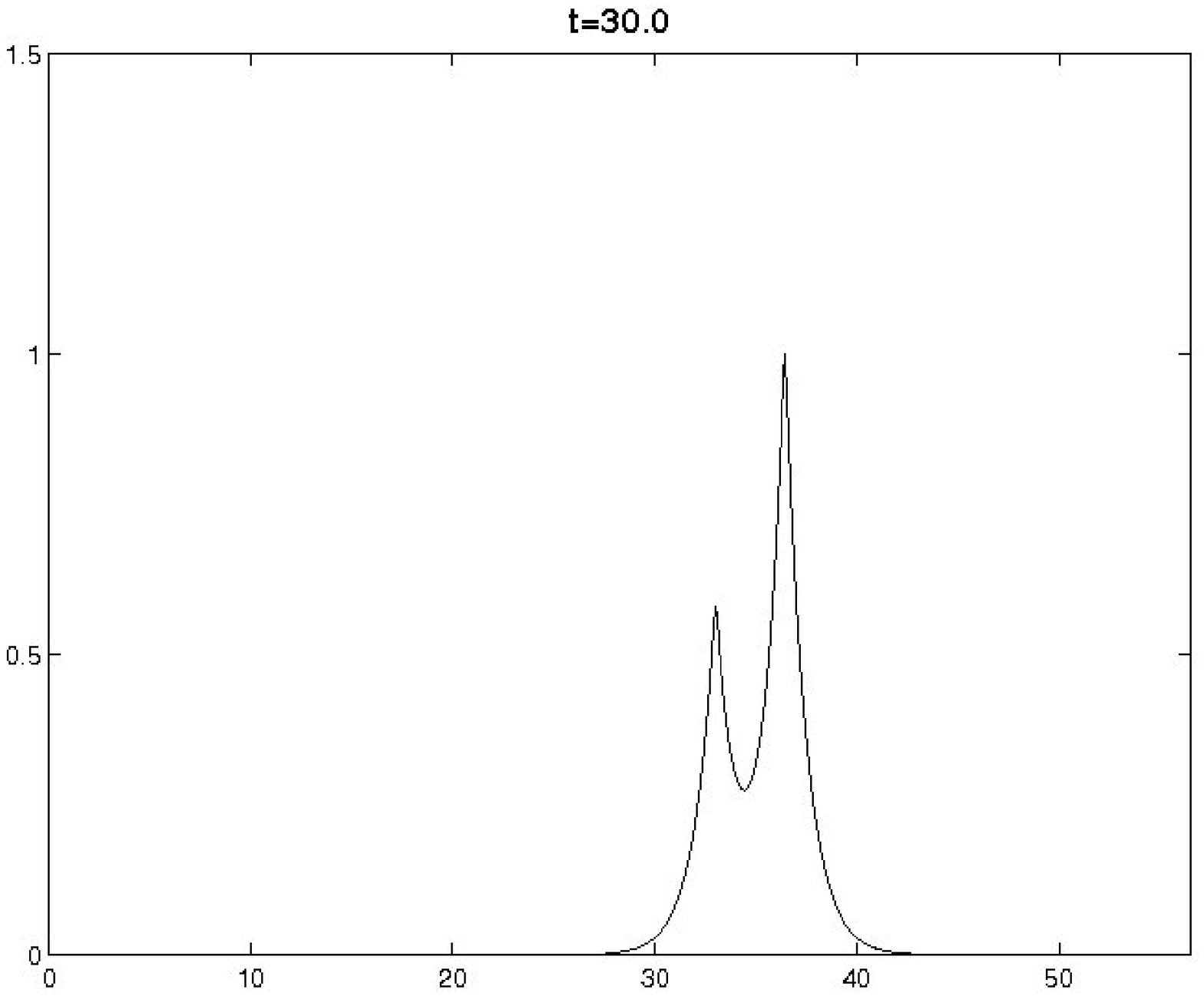}}
\scalebox{0.45}{\includegraphics{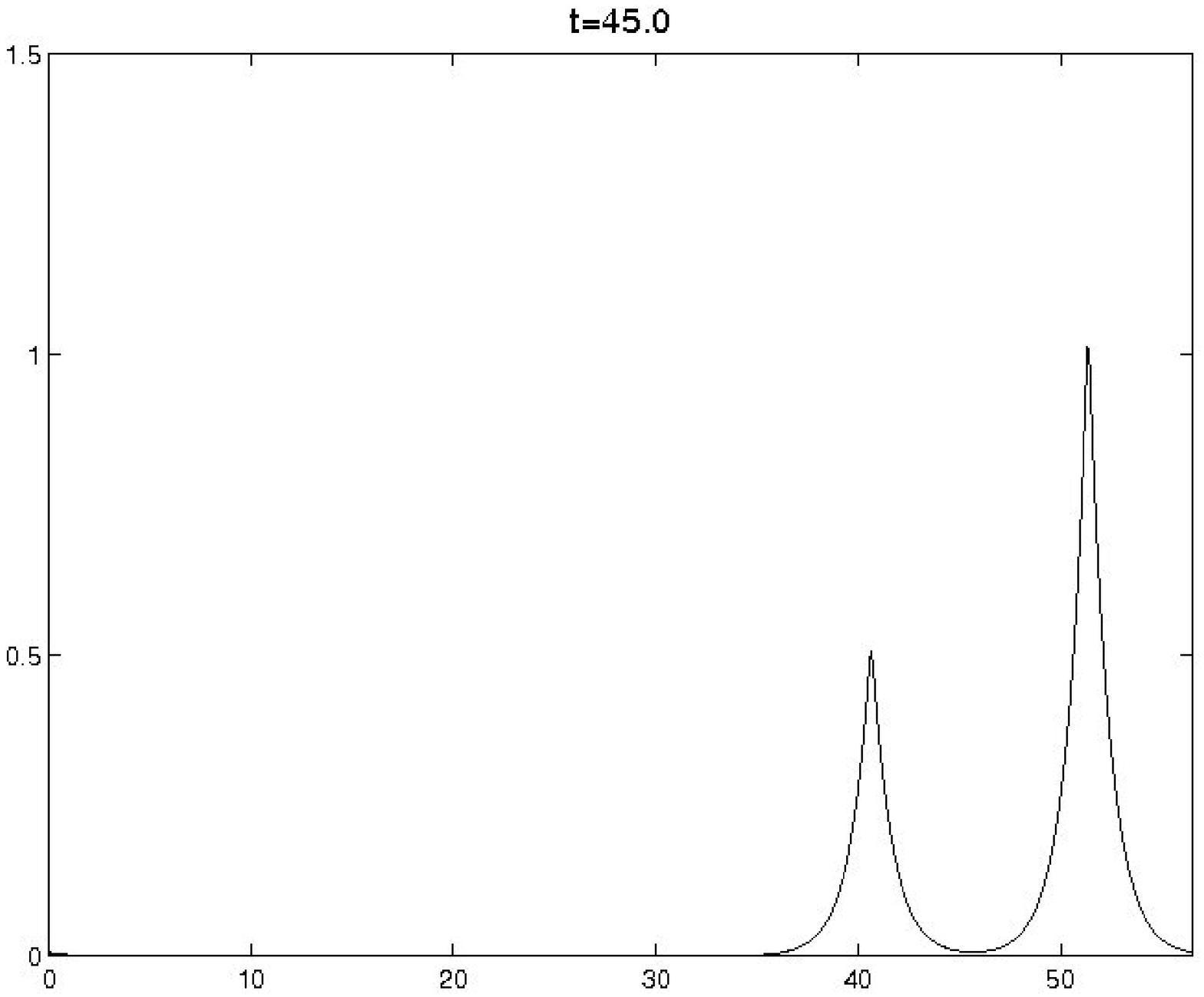}}
\end{center}
\caption{\label{overtake} Plots of the velocity during an overtaking
  collision between two peakons at times $t=1$, $t=25$, $t=35$, and
  $t=45$. The peakon to the left has greater velocity than the one to
  the right, so they collide. The momentum is transferred from the
  peakon behind to the one in front which accelerates away.}
\end{figure}

\begin{figure}[htp]
\begin{center}
\scalebox{0.7}{\includegraphics{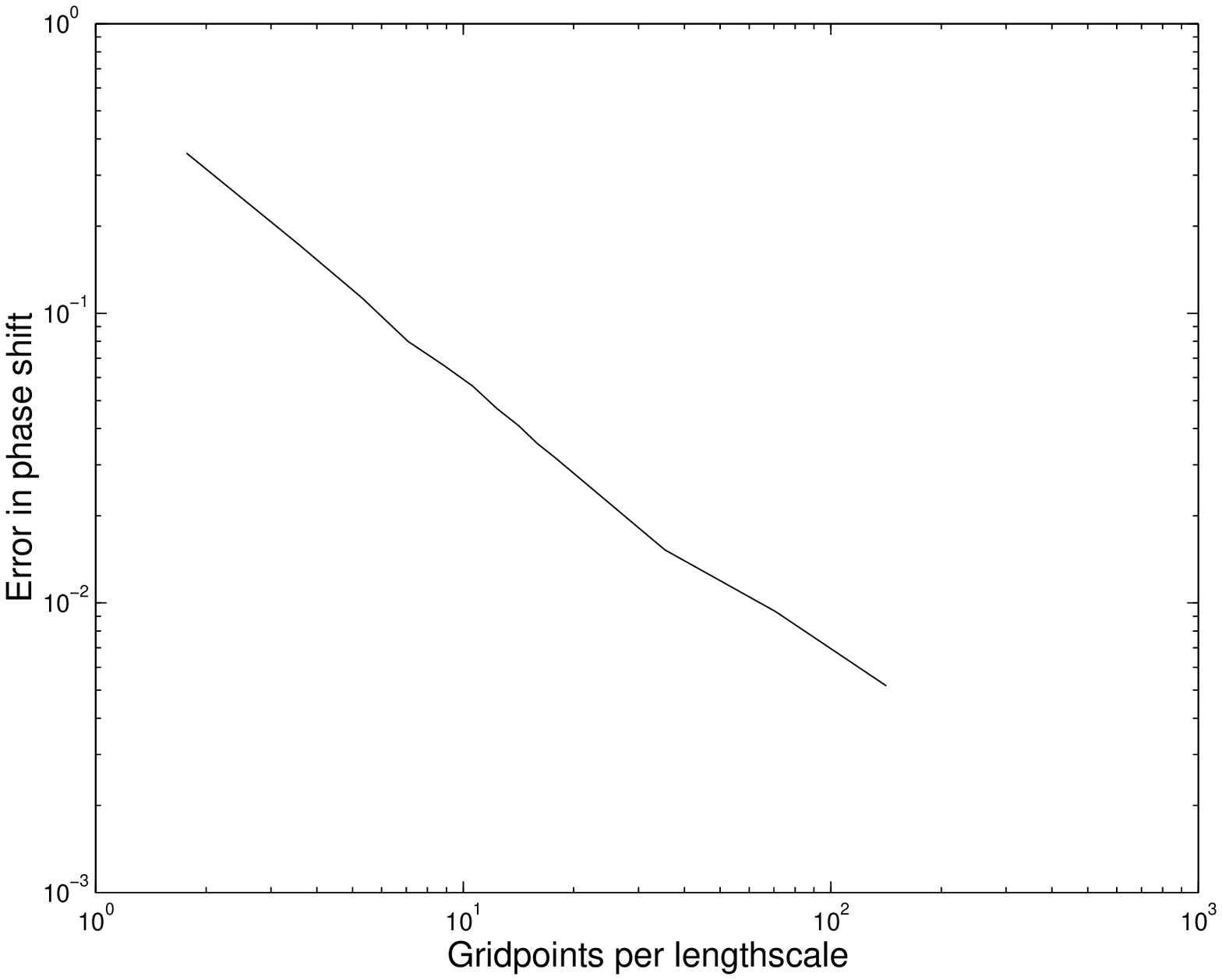}}
\end{center}
\caption{\label{phaseshift} Plot of error in calculating the fast
  shift in the faster peakon between a peakon of height 1.0 and a
  peakon of height 0.5 against grid resolution (measured as number of
  gridpoints in one characteristic length $\alpha=1$). This plot shows
  linear convergence of the error in the numerical solution.}
\end{figure}

For our second convergence test we used the problem of an overtaking
collision (illustrated in figure \ref{overtake}) between two
right-propagating peakons. \cite{CaHo1993} gives a formula the
phase-shifts for such a collision (\emph{i.e.} the asymptotic
difference in positions for the larger and smaller soliton with and
without the collision). A plot of the error in the phaseshift against
grid-size is given in figure \ref{phaseshift}.

We found that the performance of the method when solving for head-on
peakon/anti-peakon interactions was quite poor. During the collision
the two peakons approach each other and stick together once they are
both within a grid width of each other. This appears to be an issue
with representing the momentum using Lagrangian particles, as to
achieve a method with the correct results for the collision, the
particle momenta would need to go to infinity during the collision.
However, we can also view this as a benefit of the method. The
peakon/anti-peakon solution represents an instability in the
equations; namely that a small perturbation of the solution can result
in a peakon/anti-peakon pair being created. As our method does 
not support this solution at the moment of collision, this type 
of instability does not pollute our numerical results.

\subsection*{2D Flows}
In this section we show a few results obtained using the VPM
method to discretise EPDiff in two dimensions with Lagrangian
\[
L = \frac{1}{2}\int |\MM{u}|^2 + \alpha^2|\nabla\MM{u}|^2\diff{^2x},
\]
for a constant lengthscale $\alpha$ so that the velocity $\MM{u}$ is
obtained from the momentum $\MM{m}=\dd{\MM{l}}{\MM{u}}$ by inverting
the modified Helmholtz operator
\[
\MM{m} = (1-\alpha^2\Delta)\MM{u}.
\]

\begin{figure}[htp]
\begin{center}
  \scalebox{0.35}{\includegraphics{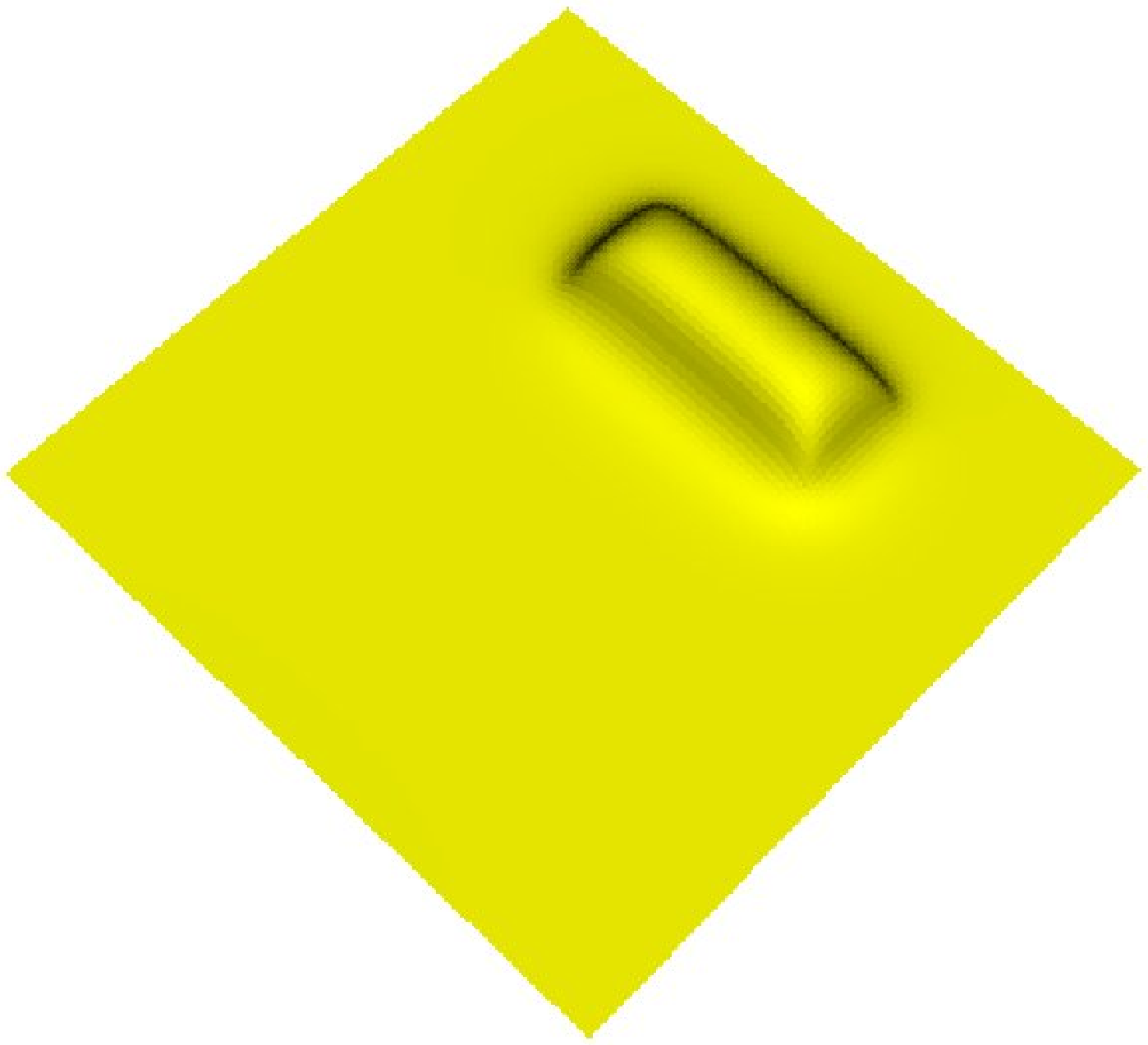}}
  \scalebox{0.35}{\includegraphics{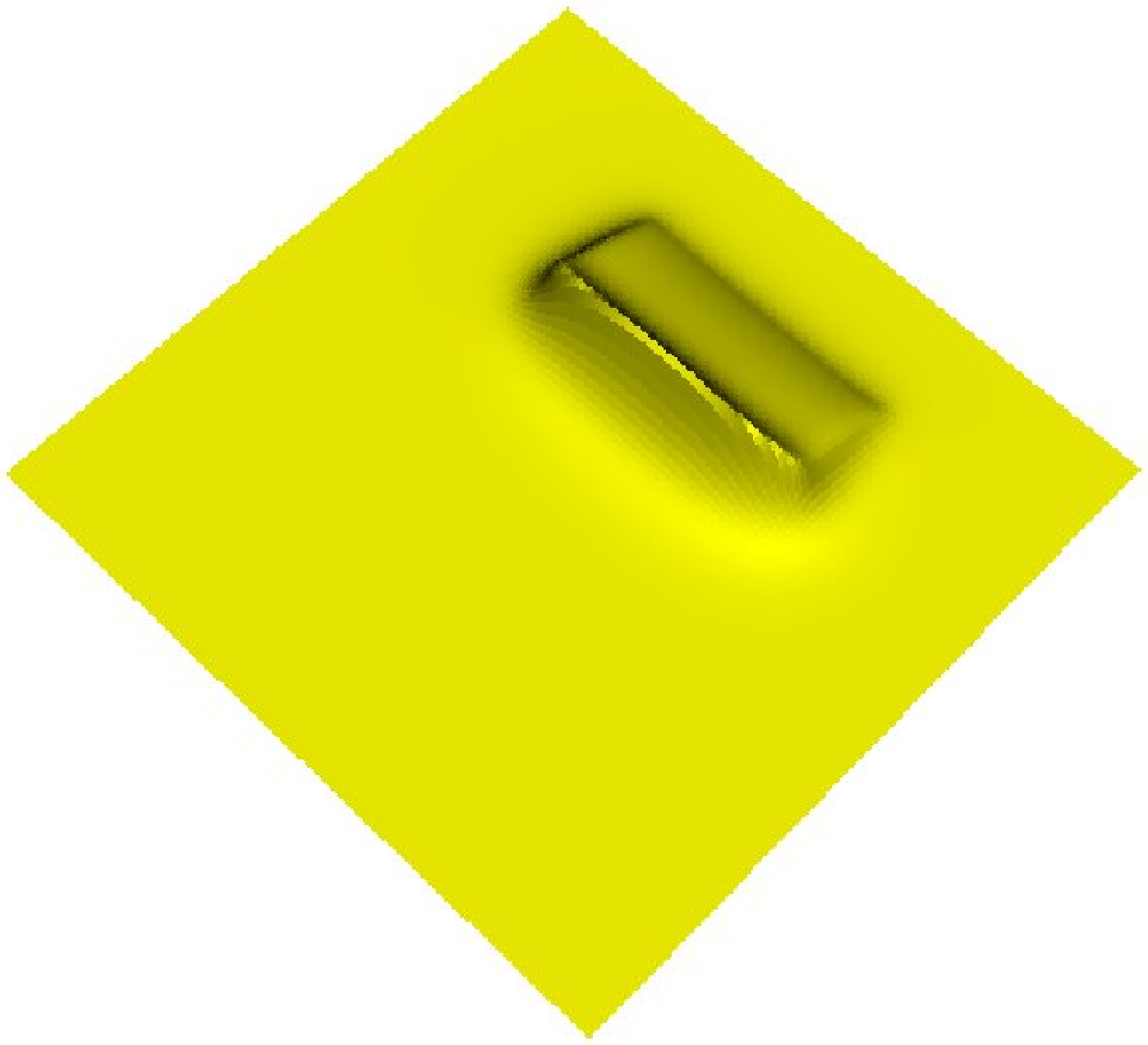}}
  \scalebox{0.35}{\includegraphics{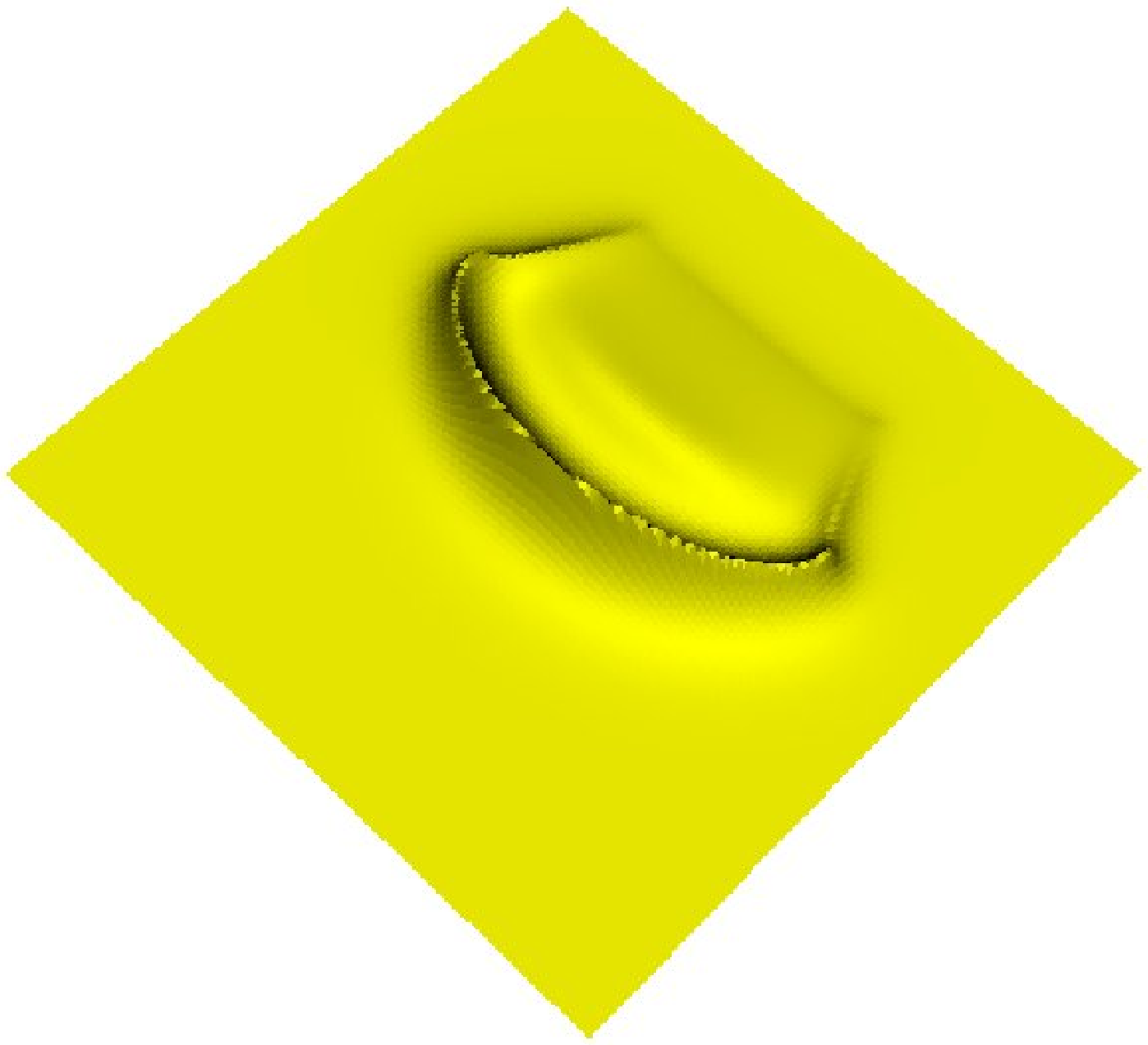}}
  \scalebox{0.35}{\includegraphics{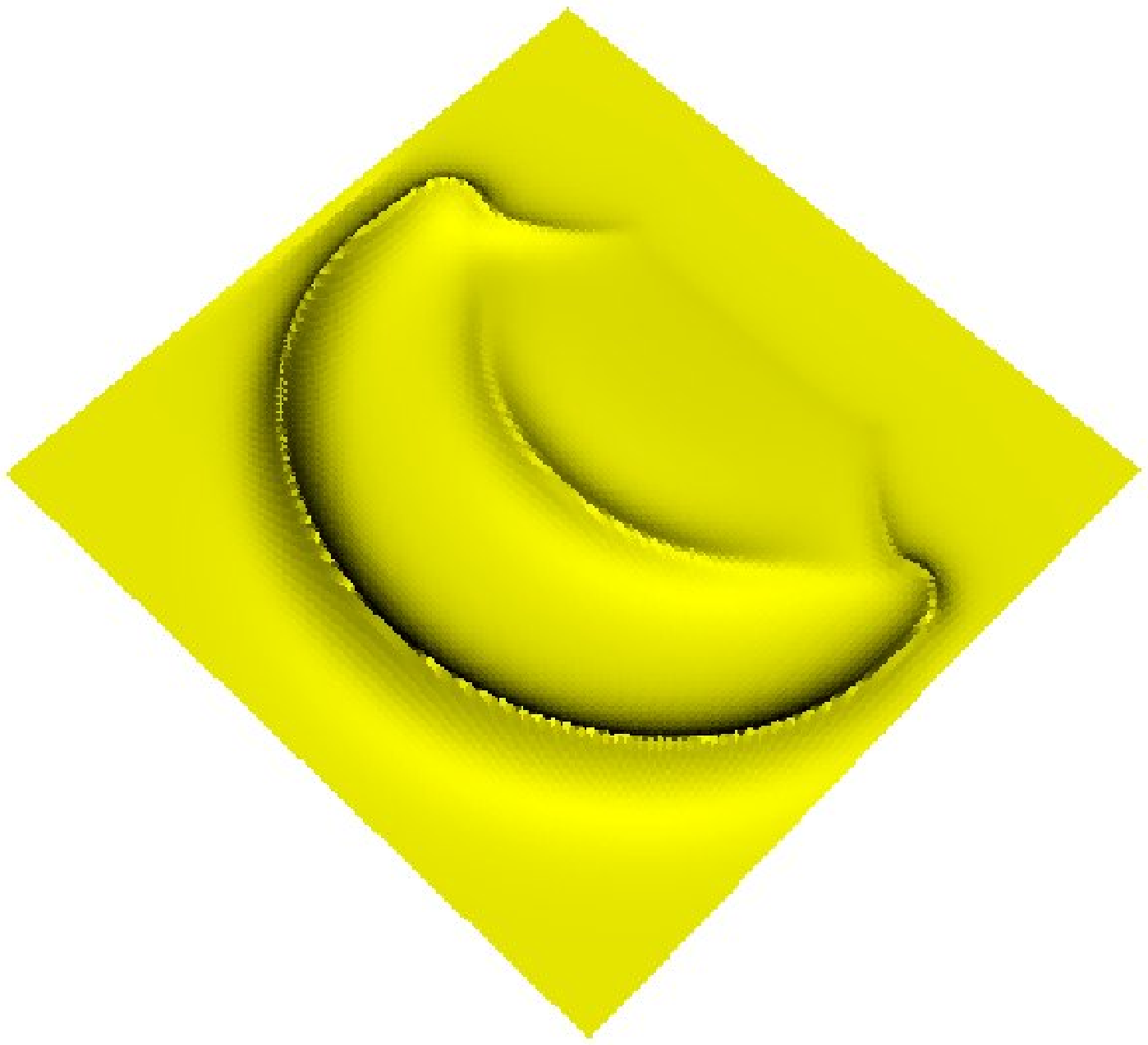}} 
\end{center}
\caption{\label{cont}Plots showing surfaces of velocity magnitude
  $|\MM{u}|$ at times $t=0$, 0.45, 1.5 and 3.4 with a $C^2$-smooth,
  compactly-supported initial condition for $\MM{u}$.  The momentum
  becomes supported on lines (so that the velocity has a ``peaked''
  profile) which spread out. This is the 2-dimensional version of
  emerging peakons illustrated in figure \ref{sech_initial_value}.
  These results were obtained using 65536 particles, a $128\times 128$
  grid with periodic boundary conditions, tensor product B-spline
  basis functions and piecewise-linear finitelement discretisation of
  the Lagrangian. The timestep is $1.0\times10^{-3}$,and $\alpha =
  0.2$.}
\end{figure}

\begin{figure}[htp]
\begin{center}
  \scalebox{0.35}{\includegraphics{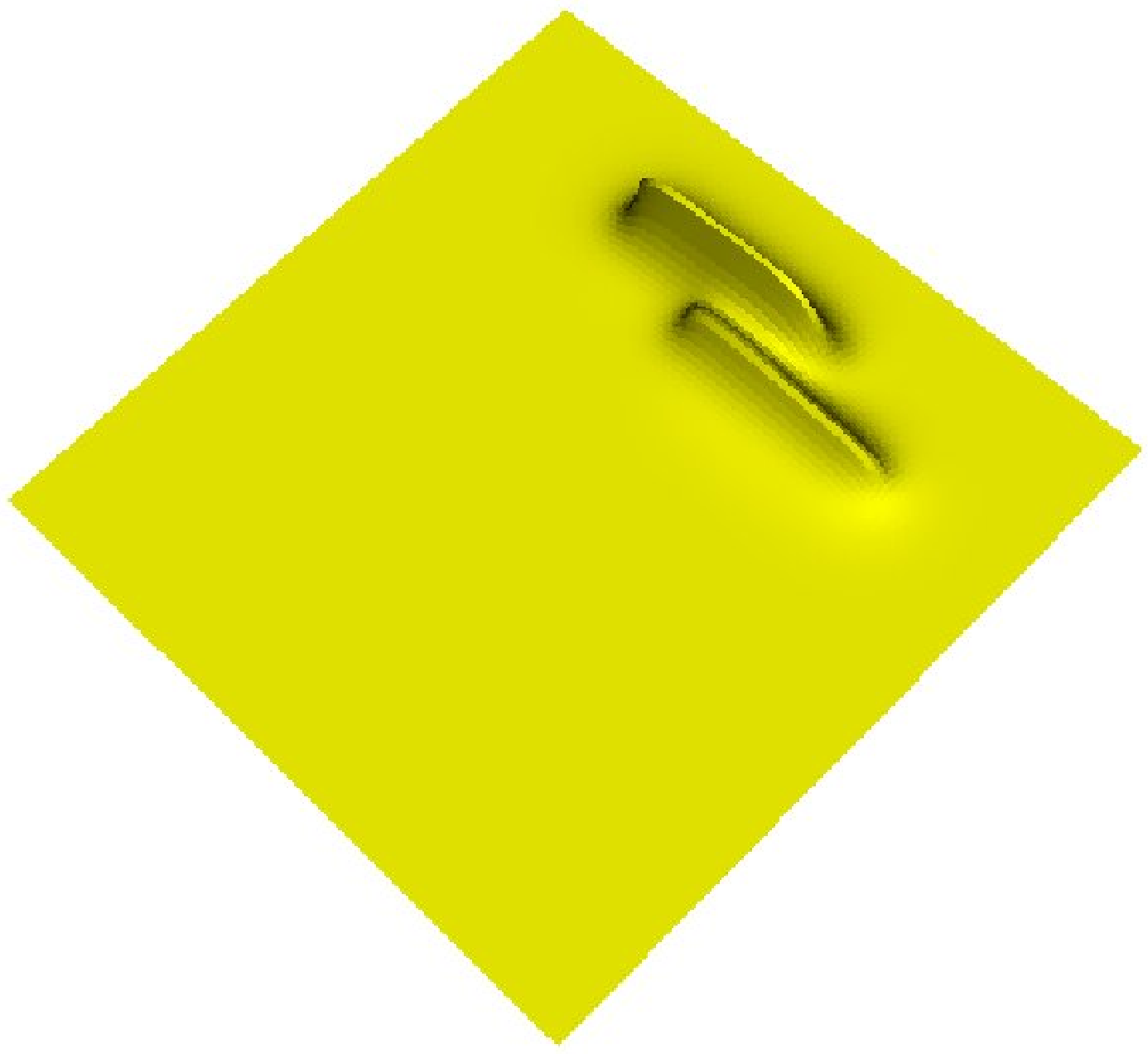}}
  \scalebox{0.35}{\includegraphics{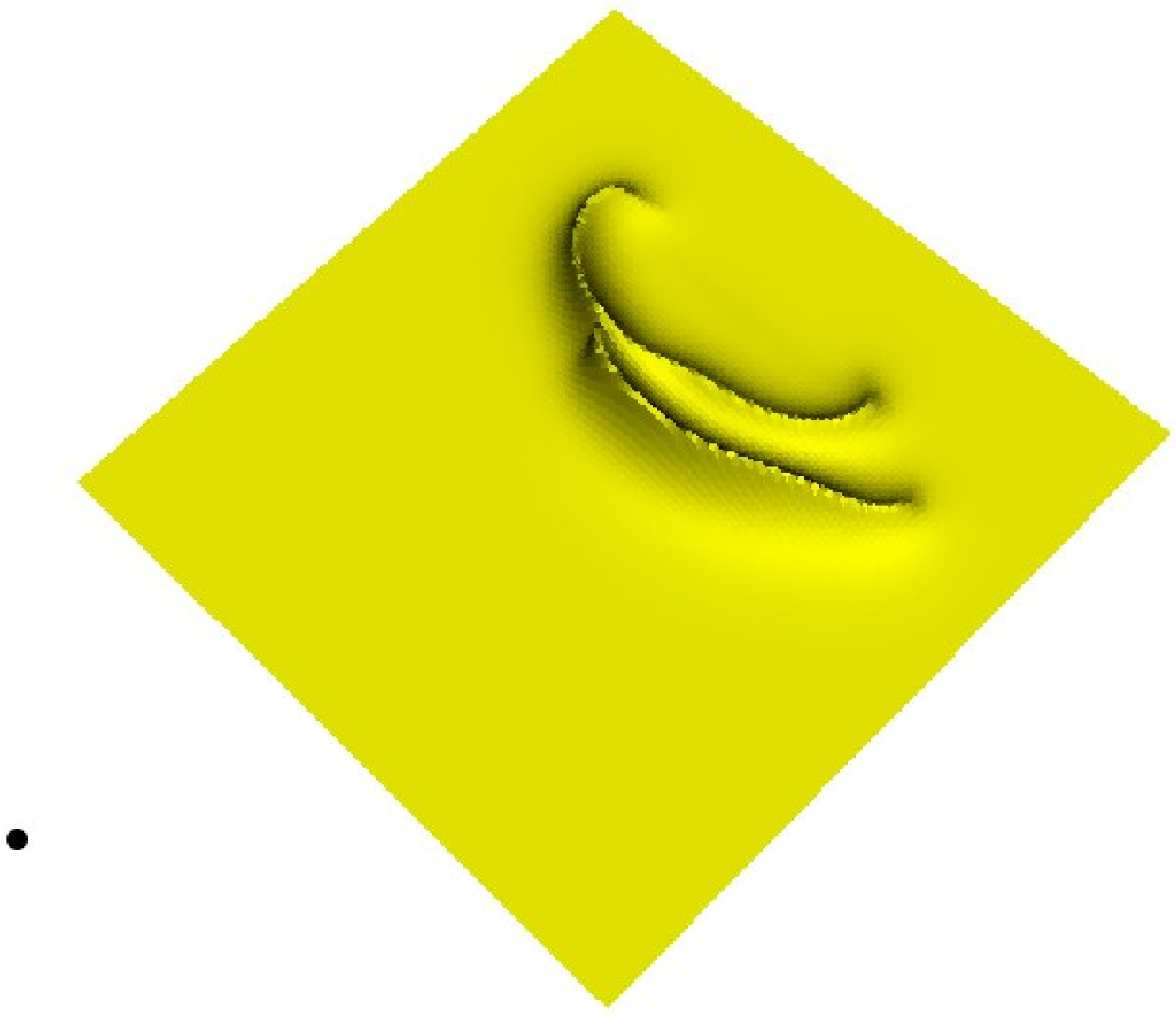}}
  \scalebox{0.35}{\includegraphics{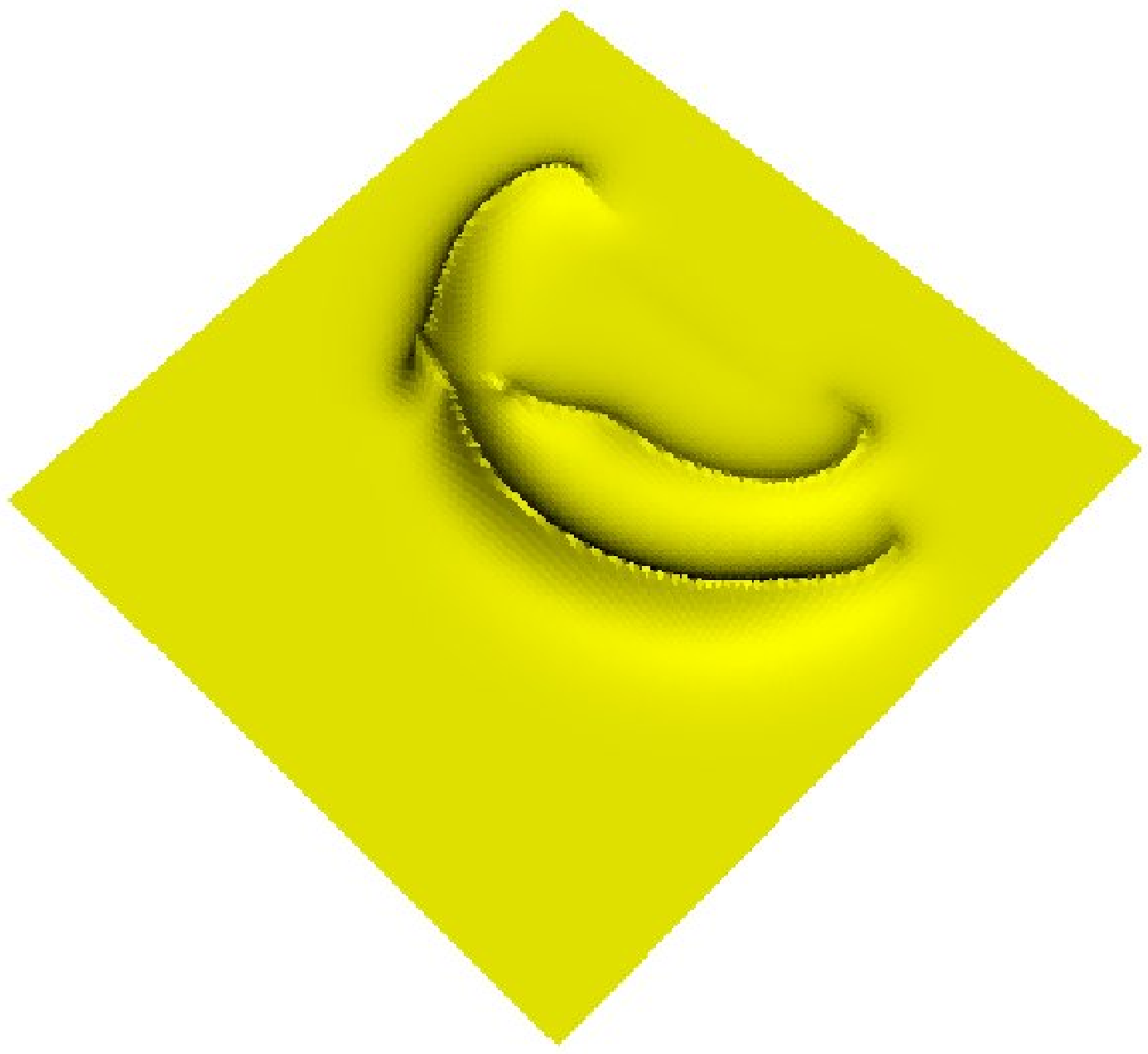}}
  \scalebox{0.35}{\includegraphics{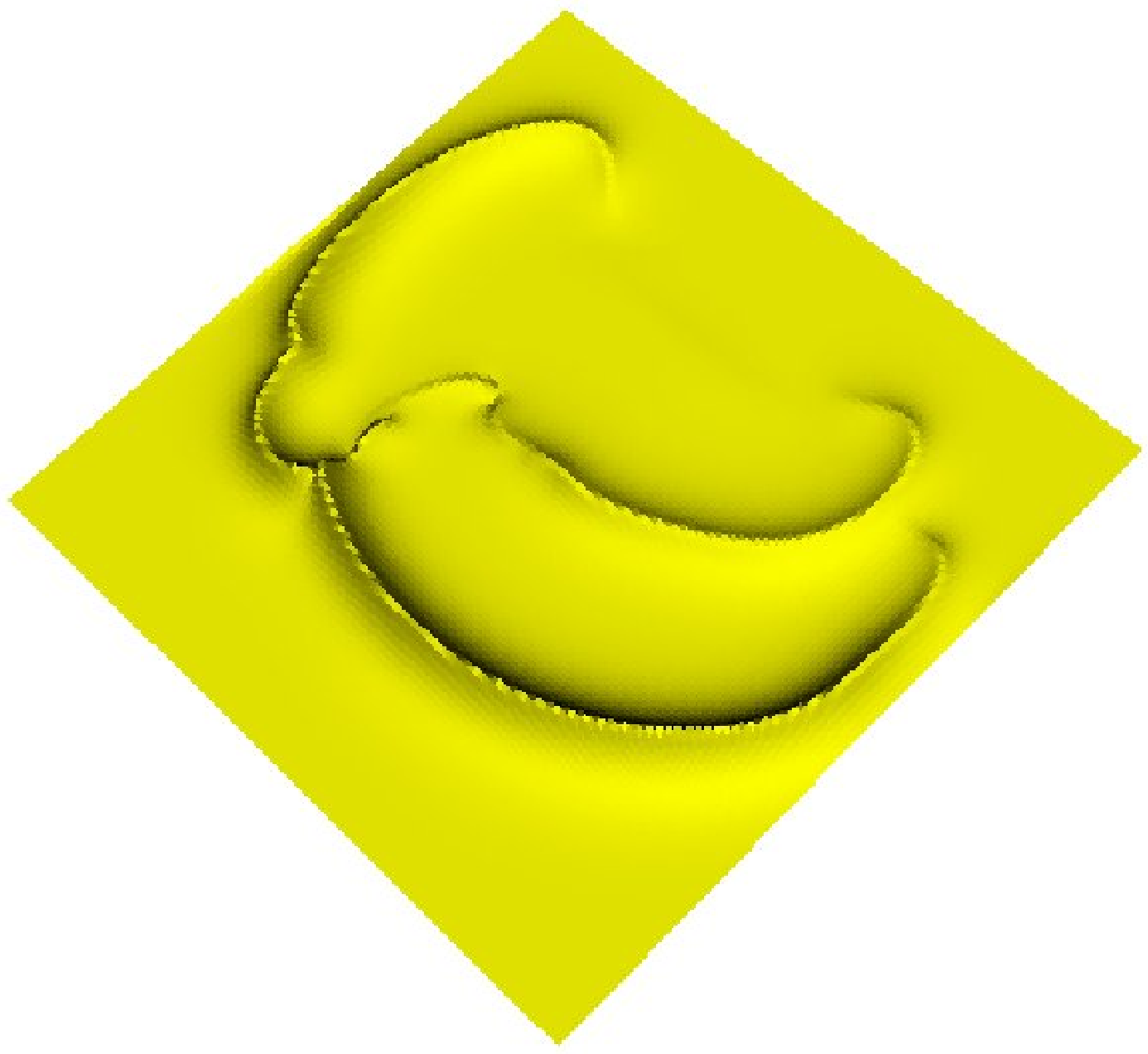}} 
\end{center}
\caption{\label{rearend}Plots showing surfaces of velocity magnitude
  $|\MM{u}|$ at times $t=0$, 1.55, 2.7 and 4.75 showing an
  ``overtaking'' collision between two singular momentum filaments.
  The filament which is initially behind has greater momentum and so
  it catches up with the filament in front, transferring momentum to
  the front filament and causing a reconnection to occur. This is the
  2-dimensional version of the process illustrated in figure
  \ref{overtake}. This is the nonlinear reconnection process which is
  illustrated in the Space Shuttle image in figures \ref{dongsha} and
  \ref{zoom dongsha}. These results were obtained using the same
  method as figure \ref{cont}.}
\end{figure}

In the first experiment the initial condition for the momentum had a
2-dimensional ``top-hat'' profile
\[
\MM{m} = (m(x,y),0), \quad m = \left\{
\begin{array}{c c}
1 & \mbox{ if } a < x < b, \quad c< y <d, \\
0 & \mbox{ otherwise.}
\end{array},
\right.
\]
for constants $a,b,c,d$, so that the velocity has continuous gradients
and has compact support. Figure \ref{cont} shows the evolution of the
velocity at subsequent times; it illustrates how EPDiff evolves to
form singular filaments of momentum from smooth initial conditions. 

\noindent {\bfseries Verification of conservation laws}
The next set of numerical results demonstrate the conservation of the
right-action momenta given in section \ref{R-Action-sec} and the
connection with Kelvin's circulation theorem. 

\begin{figure}[htp]
\begin{center}
  \scalebox{0.8}{\includegraphics{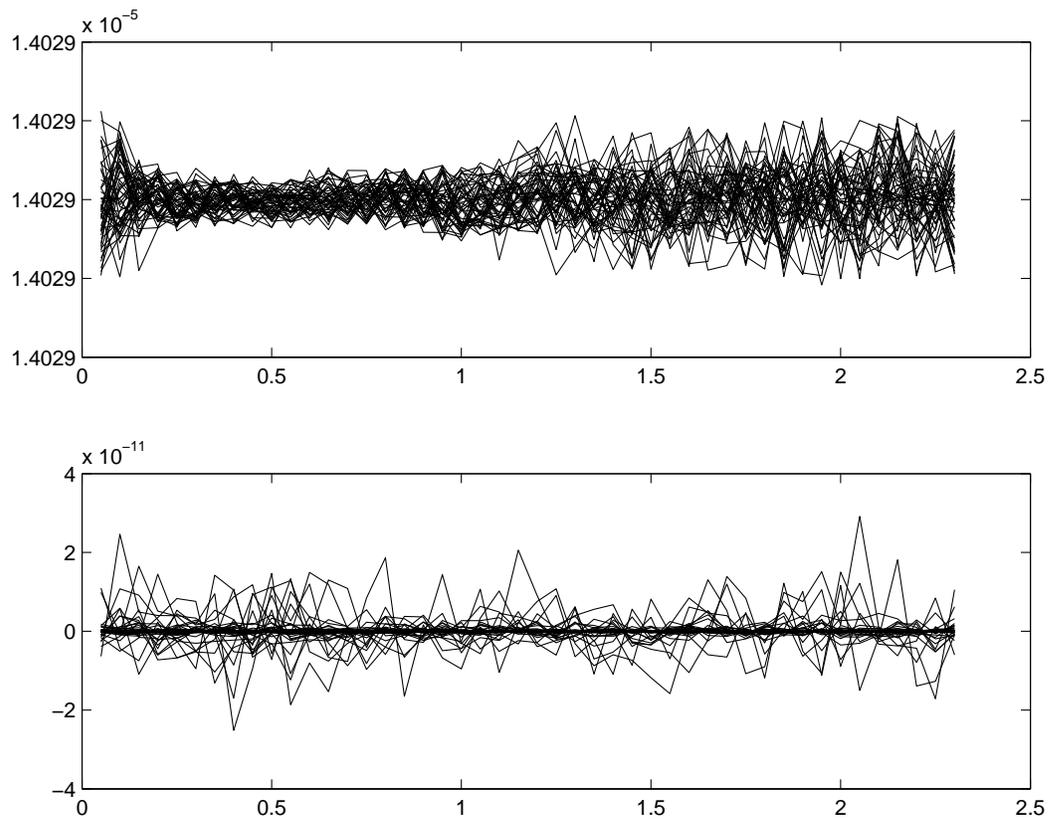}}
\end{center}
\caption{\label{PdQ} Plots of $\MM{P}_\beta\cdot J_\beta$ against time
  for a selection for particles with label $\beta$, illustrating that
  it is exactly conserved during the simulation. The upper plot is the
  $x$-component and the lower plot is the $y$-component. Any variation
  seen is due to numerical round-off error.}
\end{figure}

Figure \ref{PdQ} shows the value of the momentum map for right action
$\MM{P}_\beta\cdot J_\beta$ for a selection of particles from the flow
in figure \ref{cont}. All the particles have the same value because
they all have the same initial momentum and $J$ is set to the identity
initially. The figure shows that the numerical method preserves these
conserved momenta up to round-off error. This follows from theorem
\ref{Inv-thm} for the time continuous equations and the fact that
conserved momentum maps are also conserved by variational integrators.

To verify the discrete circulation conservation discussed in section
\ref{kelvin}, we took an arbitrary loop containing some of the
particles and advected the loop with the flow shown in figure
\ref{cont}, using
\[
\MM{Q}^{n+1}(s) = \MM{Q}^{n}(s) + \sum_k\MM{u}_k\psi_k(\MM{Q}(s)),
\]
where $s$ parameterises the loop.  During the course of the flow this
arbitrary circulation loop evolves, changing shape and length
significantly. However, the circulation around the loop remains
constant (up to numerical round off), as verified numerically in the
following.

\begin{figure}[htp]
\begin{center}
  \scalebox{0.8}{\includegraphics{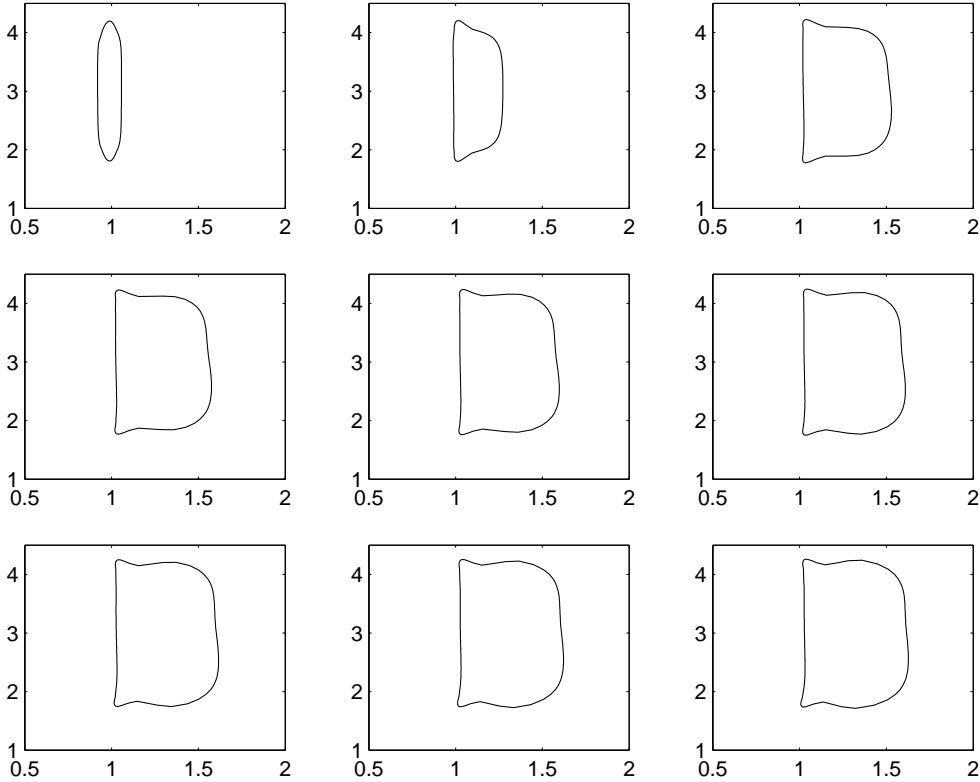}}
\end{center}
\caption{\label{loop trajs}Plots showing an embedded loop in the
  time-varying flow obtained from a solution of EPDiff, illustrated in
  figure \ref{cont}, at times (left-to-right, top-to-bottom) $t=0$,
  0.25, 0.5, 0.75, 1, 1.25, 1.5, 1.75, and 2.0.}
\end{figure}

\begin{figure}[htp]
\begin{center}
  \scalebox{0.7}{\includegraphics{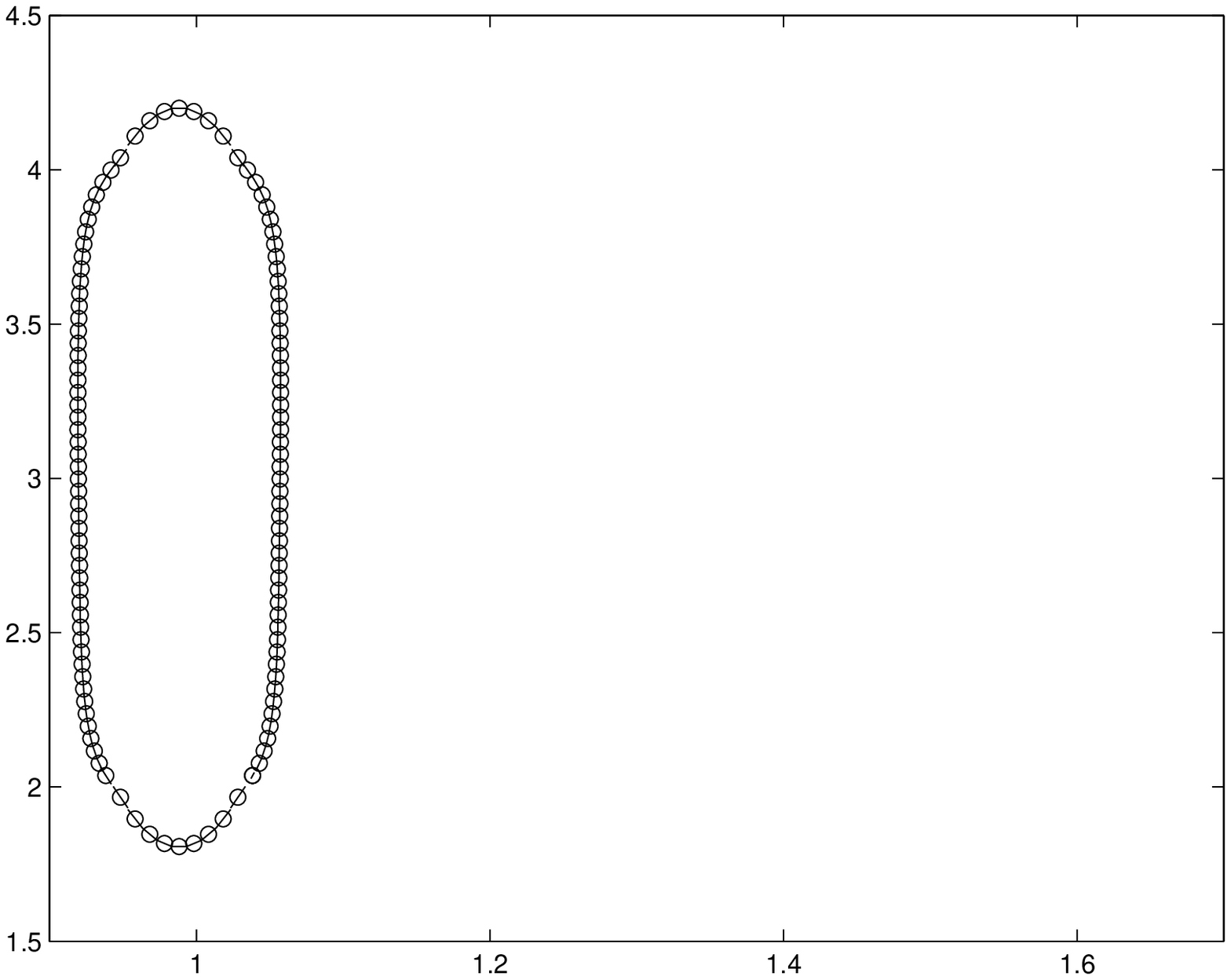}}
  \scalebox{0.7}{\includegraphics{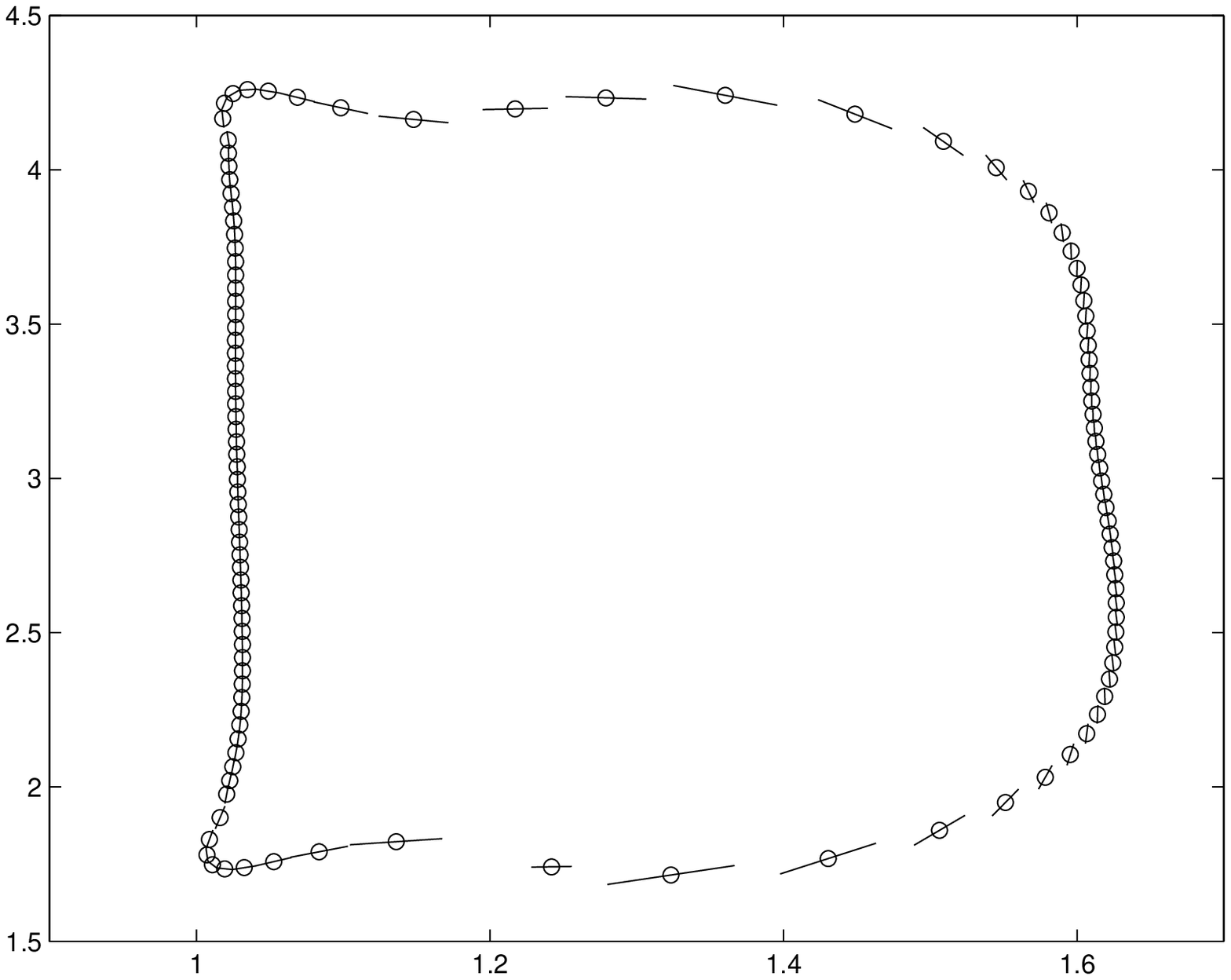}}
\end{center}
\caption{\label{looptans}Plot showing the curve with embedded line
  elements at time $t=0$ and $t=2.3$. The line elements at time
  $t=2.3$ are the exact tangents to the curve which passes through the
  tangents given at time $t=0$ and is then advected using the time
  independent flow given in figure \ref{cont}.}
\end{figure}

\begin{figure}[htp]
\begin{center}
  \scalebox{0.7}{\includegraphics{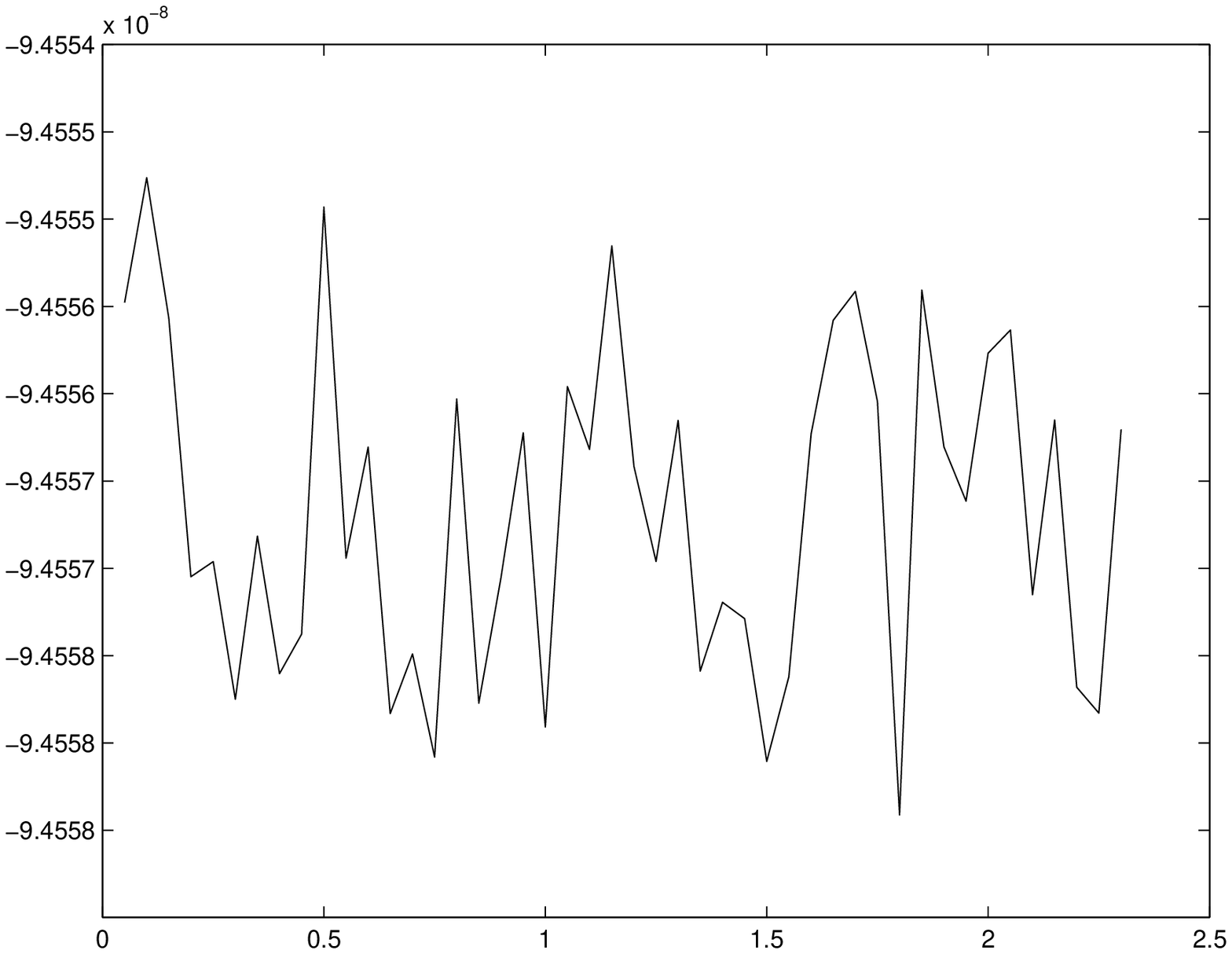}}
\end{center}
\caption{\label{circul}Plot of the circulation integral around the
  loop illustrated in figure \ref{looptans} against time. The total
  circulation is very small because the momentum is initially almost
  exactly tangent to the curve, and the conservation of
  $\MM{P}d\MM{Q}$ ensures that it stays tangent during the whole
  simulation. The circulation is exactly preserved by the numerical
  method; any variations which can be seen here are due to round-off
  error in the numerical scheme and in the calculation of the loop
  integral.}
\end{figure}

To write down the circulation integral, we choose an initial density
\[
D_k = (M^{-1})_{kl}\sum_\beta\tilde{D}_\beta\psi_k(\MM{Q}_\beta),
\]
where the values of $\tilde{D}$ do not matter much, as they are not
coupled with the dynamics. Hence, we choose the values
$\tilde{D}_\beta=1$, $\beta=1,\ldots,N_p$. To obtain the discretised loop
integral
\[
\sum_{\beta=1}^N\frac{\MM{P}_\beta}{\tilde{D}_\beta}\cdot\Delta
\MM{x}_\beta,
\]
we need to calculate $\Delta\MM{x}_\beta$ as discussed in section
\ref{kelvin}. This is done by finding discrete line elements
$\Delta\MM{x}_\beta^0$ for the initial loop and then calculating
$\Delta\MM{x}$ for subsequent timesteps using
\[
\Delta\MM{x}_\beta^n = J_\beta^n\cdot\Delta\MM{x}_\beta^0.
\]
{\bfseries Summary of circulation loop figures}
\begin{itemize}
\item Figure \ref{loop trajs} shows the
evolution in time of this circulation loop, under the flow induced by
the expanding waves in figure \ref{cont}. 
\item A plot showing initial and advected line elements is given in figure
\ref{looptans}. This plot illustrates how the line elements evolve
when a loop is stretched out by the flow. In the top and the bottom
of the loop, where the stretching is greatest, one can see how
the line elements extend to provide a numerical approximation
to $d\MM{x}$ on the loop.
\item Finally a plot of the circulation integral is given in figure
\ref{circul}. This plot shows that the circulation round the loop 
is exactly preserved during the simulation (up to round-off error
in the calculation of the discrete integral).
\end{itemize}

\section{Summary and Outlook}
In this paper we studied the Variational Particle-Mesh method
applied to the EPDiff equation. We introduced a constrained variational
principle for the method and gave discrete Euler-Poincar\'e formulae
on the Eulerian grid resulting from the variational principle which
show that the grid velocities and momenta satisfy the EPDiff equations
in Eulerian form. Next we looked at left- and right-actions of
velocity vector fields on the Lagrangian particles and obtained
corresponding momentum maps. The left-action, when restricted to the
finite space of velocity fields used in the method, gives rise to a
momentum map which is the same formula as used for calculating the
grid momentum from the particle variables. The right-action, which had
to be interpreted in a wider space, can be interpreted as a discrete
form of particle-relabelling since it corresponds to moving the
particles and also changing the momenta in such a way so that the grid
velocities remain constant. Finally we gave some interpretation of
these transformations in terms of matrices which determine the
local deformation of infinitesimal line elements, thereby allowing us to 
write down discrete loop integrals on advected loops. This led to
a discrete circulation theorem.

Our next aim is to find an extension of this work which gives a
discrete circulation theorem for fluid PDEs which involve mass density
and other advected quantities as well as velocity. The general
approach, following the continuous theory, will be to
\begin{itemize}
\item specify the transformations corresponding to discrete
  relabelling,
\item determine the transformation of density and other advected
  quantities under this discrete relabeling group,
\item calculate the momentum densities obtained from these transformations.
\item show that the ratio of these momentum densities to the mass
  density is invariant.
\end{itemize}
Including advected quantities in this awy will allow introduction of
potential energy and hence linear dispersion effects into the
numerical description of the internal wave interactions using the VPM
method. 

\subsection*{Acknowledgements} We are grateful to our colleagues Joel
Fine and Matthew Dixon at Imperial College London for their advice and
consultation regarding this problem. We are also grateful to the
ONR-NLIWI program for partial funding of this endeavor, and to Tony
Liu for use of the SAR images of the South China Sea taken from the
Space Shuttle. DDH is also grateful for partial support from the
Office of Science, US Department of Energy.


\nocite{*}

\bibliography{LatticeEPDiff1}

\end{document}